\documentclass[12pt]{article}
\usepackage{amssymb,amsmath,bm,graphicx}
\usepackage{mathrsfs}
\usepackage{constants}
\usepackage{hyperref}
\usepackage{color}

\date{}
\begin{document}
\newcommand{\bea}{\begin{eqnarray}}
\newcommand{\ena}{\end{eqnarray}}
\newcommand{\beas}{\begin{eqnarray*}}
\newcommand{\enas}{\end{eqnarray*}}
\newcommand{\beq}{\begin{equation}}
\newcommand{\enq}{\end{equation}}
\def\qed{\hfill \mbox{\rule{0.5em}{0.5em}}}
\newcommand{\bbox}{\hfill $\Box$}
\newcommand{\ignore}[1]{}
\newcommand{\ignorex}[1]{#1}
\newcommand{\wtilde}[1]{\widetilde{#1}}
\newcommand{\qmq}[1]{\quad\mbox{#1}\quad}
\newcommand{\qm}[1]{\quad\mbox{#1}}
\newcommand{\nn}{\nonumber}
\newcommand{\Bvert}{\left\vert\vphantom{\frac{1}{1}}\right.}
\newcommand{\To}{\rightarrow}
\newcommand{\E}{\mathbb{E}}
\newcommand{\Var}{\mathrm{Var}}
\newcommand{\Cov}{\mathrm{Cov}}
\newcommand{\Corr}{\mathrm{Corr}}
\makeatletter
\newsavebox\myboxA
\newsavebox\myboxB
\newlength\mylenA
\newcommand*\xoverline[2][0.70]{%
    \sbox{\myboxA}{$\m@th#2$}%
    \setbox\myboxB\null
    \ht\myboxB=\ht\myboxA%
    \dp\myboxB=\dp\myboxA%
    \wd\myboxB=#1\wd\myboxA
    \sbox\myboxB{$\m@th\overline{\copy\myboxB}$}
    \setlength\mylenA{\the\wd\myboxA}
    \addtolength\mylenA{-\the\wd\myboxB}%
    \ifdim\wd\myboxB<\wd\myboxA%
       \rlap{\hskip 0.5\mylenA\usebox\myboxB}{\usebox\myboxA}%
    \else
        \hskip -0.5\mylenA\rlap{\usebox\myboxA}{\hskip 0.5\mylenA\usebox\myboxB}%
    \fi}
\makeatother

\newtheorem{theorem}{Theorem}[section]
\newtheorem{corollary}[theorem]{Corollary}
\newtheorem{conjecture}[theorem]{Conjecture}
\newtheorem{proposition}[theorem]{Proposition}
\newtheorem{lemma}[theorem]{Lemma}
\newtheorem{definition}[theorem]{Definition}
\newtheorem{example}[theorem]{Example}
\newtheorem{remark}[theorem]{Remark}
\newtheorem{case}{Case}[section]
\newtheorem{condition}{Condition}[section]
\newcommand{\proof}{\noindent {\it Proof:} }

\title{{\bf\Large Stein's method using approximate zero bias couplings with applications to combinatorial central limit theorems under the Ewens distribution}}
\author{Nathakhun Wiroonsri \thanks{Department of Mathematics, University of Southern California, Los Angeles, CA 90089, USA. Email: wiroonsr@usc.edu }  \\ University of Southern California}

\footnotetext{Research supported by the 2017-2018 Russell Endowed Fellowship, University of Southern California.}
\maketitle

\begin{abstract} 
We generalize the well-known zero bias distribution and the $\lambda$-Stein pair to an approximate zero bias distribution and an approximate $\lambda,R$-Stein pair, respectively. Berry Esseen type bounds to the normal, based on approximate zero bias couplings and approximate $\lambda,R$-Stein pairs, are obtained using Stein's method. The bounds are then applied to combinatorial central limit theorems where the random permutation has the Ewens $\mathcal{E}_\theta$ distribution with $\theta>0$ which can be specialized to the uniform distribution by letting $\theta=1$. The family of the Ewens distributions appears in the context of population genetics in biology. 
\end{abstract}

\section{Introduction} \label{Intro}
We develop $L^1$ and $L^\infty$ bounds for normal approximation using Stein's method, based on approximate zero bias couplings. The results are applied to the combinatorial central limit theorems, that is, we derive such bounds for the distribution, introduced in \cite{Hoe51}, of 
\bea \label{combdef2}
Y =\sum_{i=1}^n a_{i,\pi(i)}
\ena
where $A \in \mathbb{R}^{n \times n}$ is a given real matrix with components $\{a_{i,j} \}_{i,j=1}^n$ and $\pi \in \mathcal{S}_n$ has the Ewens distribution. We recall that the $L^1$ and $L^\infty$ distances between the distributions ${\cal L}(X)$ and ${\cal L}(Y)$ of real valued random variables $X$ and $Y$ are given, respectively, by
\bea \label{Wasdef}
d_1 \big({\cal L}(X),{\cal L}(Y)\big) &=& \int_{-\infty}^\infty |P(X \le t) - P(Y \le t)| dt \nn \\
                                      &=& \sup_{h \in \mathcal{H}_1} |\E h(X)-\E h(Y)| 
\ena 
where $\mathcal{H}_1 =\{h: |h(y)-h(x)| \le |y-x|\}$, and 
\bea \label{Koldef}
d_\infty \big({\cal L}(X),{\cal L}(Y)\big) &=& \sup_{\-\infty < t < \infty} |P(X<t) - P(Y<t)| \nn \\
                                           &=& \sup_{h \in \mathcal{H}_\infty} |\E h(X)-\E h(Y)| 
\ena 
where $\mathcal{H}_\infty =\{\mathbf{1}(\cdot \le t): t \in \mathbb{R}\}$. In the following, we will drop the subscripts $1$ and $\infty$ when the statement is true for both $\mathcal{H}_1$ and $\mathcal{H}_\infty$.

Stein's method for normal approximation, introduced by Charles Stein in \cite{Stein72} (see also the text \cite{CGS11} and the introductory notes \cite{Ross11}), was motivated from the fact that $W$ has the standard normal distribution, denoted $\mathcal{N}(0,1)$, if and only if 
\beas
\E W f(W) = \E f'(W)
\enas
for all absolutely continuous functions $f$ with $\E |f'(W)| < \infty$. This equation and the form of the distances in \eqref{Wasdef} and \eqref{Koldef} lead to the differential equation
\bea \label{steineq}
h(w) - N h = f'_h(w)-wf_h(w)
\ena
where $Nh = \E h(Z)$ with $Z \sim \mathcal{N}(0,1)$ and $h \in \mathcal{H}$. Taking the supremum over all $h \in \mathcal{H}_1$ (resp. $h \in \mathcal{H}_\infty$) to the expectation on the left hand side of \eqref{steineq} with $w$ replaced by a variable $W$ yields the distance between $W$ and $Z$ in \eqref{Wasdef} (resp. \eqref{Koldef}). Thus, instead of working on the distances directly, one can handle the expectation on the right hand side using the bounded solution $f_h$ of \eqref{steineq} for the given $h$. Using this device, Stein's method has uncovered an alternative way to show convergence in distribution with additional information on the finite sample distance between distributions and can also deal with various kinds of dependence through the help of coupling constructions. 

One of the well-known couplings in the literature is the \textit{zero bias coupling} which was first introduced in \cite{GR97}. Recall that for $X$ with mean zero and variance $\sigma^2 \in (0,\infty)$, we say that $X^z$ has the $X$-zero biased distribution if 
\bea \label{zerobias}
\E Xf(X) = \sigma^2 \E f'(X^z)
\ena
for all absolutely continuous functions $f$ for which the expectations exist. Applying \eqref{zerobias} in \eqref{steineq} with $X$ replaced by $X/\sigma$, we have 
\bea \label{genzero}
\E h(X/\sigma) - Nh &=& \E \left(f'\left(\frac{X}{\sigma}\right)-\frac{X}{\sigma} f\left(\frac{X}{\sigma} \right) \right) \nn\\
                    &=& \E \left(f'\left(\frac{X}{\sigma}\right)-f'\left(\frac{X^z}{\sigma} \right) \right) . 
\ena
Generalizing the proofs of Theorems 4.1 and 5.1 of \cite{CGS11} that only showed the results for $\sigma^2=1$, with the help of \eqref{genzero}, we have 
\bea \label{oldl1}
d_1(\mathcal{L}(X/\sigma),\mathcal{L}(Z)) \le \frac{2}{\sigma}\E |X^z-X|,
\ena
and, with $|X^z-X|\le \delta$,
\bea \label{oldKol}
d_\infty(\mathcal{L}(X/\sigma),\mathcal{L}(Z)) \le \frac{(1+1/\sqrt{2\pi}+\sqrt{2\pi}/4)\delta}{\sigma}.
\ena
It is easy to see from \eqref{oldl1} and \eqref{oldKol} that once the zero bias coupling has been constructed in such a way that the two variables are close, one can simply obtain good bounds for $L^1$ and $L^\infty$ distances. Nevertheless, the difficult part is that there is no general way to construct the zero bias coupling. One of the most efficient method, introduced in \cite{GR97}, is to take advantage of the existence of a $\lambda$-Stein pair. We recall that an exchangeable pair $X',X''$ forms a $\lambda$-Stein pair if 
\beas
\E(X''|X') = (1-\lambda)X'
\enas
for some $0<\lambda<1$. The following lemma illustrates the way to construct zero bias couplings through $\lambda$-Stein pairs.
\begin{lemma}[\cite{GR97}] \label{zerobiaslemma}  
Let $X',X''$ be a $\lambda$-Stein pair with $\Var(X') = \sigma^2 \in (0,\infty)$ and distribution $F(x',x'')$. Then when $X^\dagger,X^\ddagger$ have distribution
\beas
d F^\dagger(x',x'') = \frac{(x'-x'')^2}{\E(X'-X'')^2}dF(x',x''),
\enas
and $U$ be uniform $\mathcal{U}[0,1]$ and is independent of $X^\dagger,X^\ddagger$, the variable
\beas
X^z = UX^\dagger+(1-U)X^\ddagger \text{ \ has the $X'$-zero biased distribution.}
\enas
\end{lemma}

Lemma \ref{zerobiaslemma} has been used in several works. Berry Esseen type bounds to the normal, based on zero bias couplings, were first obtained in \cite{Gol05}. The bounds were then applied to the combinatorial central limit theorems where the random permutation has either the uniform distribution or one which is constant over permutations with the same cycle type and having no fixed points. Concentration inequalities on the same setting were shown in \cite{GI14}. In \cite{FG11}, the zero biasing and this lemma were also used to obtain $L^1$ bounds to the normal for a variable constructed from an interesting property of the $\text{Jack}_\alpha$ measure on the set of partitions of size $n$. Apart from the use of this lemma, Stein's method with different techniques has been applied to the combinatorial central limit theorems under the uniform distribution in several works. One of the most recent paper is \cite{CF15} where the exchangeable pairs technique was used to obtain $L^\infty$ bounds between $Y$ as in \eqref{combdef2} with a fixed matrix $A$ replaced by independent random variables $\{X_{i,j}:i,j\in [n]\}$ and the normal distribution. Here for a positive integer $m$ we denote $[m] = \{1,2,\ldots,m \}$. The bounds there are given in term of the third moments of $X_{i,j}$. To learn more about the history of Stein's method and the combinatorial central limit theorems, see the references therein. Without Stein's method, the results were generalized to the case without third moments in \cite{Fro14} and to various moment conditions in \cite{Fro17}. The original idea of replacing $a_{i,j}$ by $X_{i,j}$ dates back to \cite{HC78} where the results are optimal only in the case that there exists $C>0$ such that $|X_{i,j}| \le C$ for all $i,j \in [n]$.

Although the combinatorial central limit theorems have attracted attention for quite some time and have been extended to different settings, the case where the random permutation has the Ewens distribution has yet been studied. It is interesting to investigate the robustness and the sensitivity of normality when the usual assumptions in the uniform distribution or the distribution that is constant over permutations with the same cycle type and having no fixed points are not satisfied. This is useful in the real-life situation as the uniform properties are sometimes believed to hold but actually do not. One difficulty that may arise in order to use Lemma \ref{zerobiaslemma} is that a $\lambda$-Stein pair does not always exist. However, one might be able to construct a pair which has nearly the same condition as the $\lambda$-Stein pair and this is where we start. We call a pair of random variables $(Y',Y'')$, an \textit{approximate $\lambda,R$-Stein pair} if it is exchangeable and satisfies
\beas
\E Y' = 0, \quad \Var(Y') = \sigma^2 
\enas
with $\sigma^2 \in (0,\infty)$, and
\bea \label{approxsteinpair}
\E(Y''|Y') = (1-\lambda)Y'+R,
\ena
for some  $0< \lambda < 1$ and $R=R(Y')$.  

In this work, we generalize Lemma \ref{zerobiaslemma} to a new version, Lemma \ref{Wdagger}, replacing a $\lambda$-Stein pair by an approximate $\lambda,R$-Stein pair. The lemma leads to a variable that is similar to the zero bias variable but has two extra terms depending on $R$ and we call it an \textit{approximate zero bias variable}. Then we also generalize the $L^1$ and $L^\infty$ bounds in \eqref{oldl1} and \eqref{oldKol} using approximate zero bias couplings in Theorems \ref{thm:approx.zb.dist.l1} and \ref{thm:approx.zb.dist.Kol}, respectively.

The remainder of this work is organized as follows. In Section \ref{main} we state and prove the general results for approximate $\lambda,R$-Stein pairs and approximate zero bias variables. These results are then applied, in Section \ref{comb}, to the combinatorial central limit theorems where the random permutation has the Ewens distribution. The description of the Ewens distribution and some necessary properties are presented in Section \ref{Ewens}. In Section \ref{appendix}, the Appendix, we prove some of the results from Section \ref{comb} that are straightforward but requires some attention to detail.


\section{Main results: approximate $\lambda,R$-Stein pairs and approximate zero bias couplings} \label{main}

Let $(Y',Y'')$ be an approximate $\lambda,R$-Stein pair.
Taking expectation in \eqref{approxsteinpair}, using exchangeability and that $Y'$ has mean zero yields
\beas
\E R=(1-\lambda)\E Y'+\E R=\E Y''=0.
\enas
In addition, for any function $f$ such that the following expectations exist,
\bea \label{approxsteinpair2}
\E Y''f(Y') &=& \E(\E(Y''f(Y')|Y')) = \E(f(Y')\E(Y''|Y')) \nn \\ 
&=& \E(f(Y')((1-\lambda)Y'+R)) = (1-\lambda)\E Y'f(Y') + \E Rf(Y').
\ena
In particular, specializing \eqref{approxsteinpair2} to the case $f(y) = y$ yields
\beas
\E Y''Y' = (1-\lambda)\E Y'^2 + \E Y'R,
\enas
and thus
\bea \label{diff2}
\E(Y''-Y')^2 &=& 2\left(\E Y'^2 - \E Y''Y' \right)\nn\\ &=& 2\left(\lambda\E Y'^2 - \E Y'R \right) = 2\left(\lambda \sigma^2 -\E Y'R \right).
\ena

Now we state and prove the following lemma which is the generalized version of Lemma \ref{zerobiaslemma} adapted to approximate $\lambda,R$-Stein pairs. We call a variable $Y^*$ that satisfies \eqref{approxzerobias} below an \textit{approximate $Y'$-zero bias variable}. 

\begin{lemma} \label{Wdagger}
Let $(Y',Y'')$ be an approximate $\lambda,R$-Stein pair with distribution $F(y',y'')$. Then when $(Y^\dagger,Y^\ddagger)$ has distribution 
\beas
dF^\dagger(y',y'') = \frac{(y''-y')^2}{\E (Y''-Y')^2}dF(y',y''),
\enas 
and $U \sim \mathcal{U}([0,1])$ is independent of $Y^\dagger,Y^\ddagger$, the variable $Y^* = UY^\dagger + (1-U)Y^\ddagger$ satisfies
\bea \label{approxzerobias}
\E [Y'f(Y')] = \sigma^2\E f'(Y^*) - \frac{\E Y'R}{\lambda}\E f'(Y^*) + \frac{\E Rf(Y')}{\lambda}
\ena
for all absolutely continuous functions $f$.
\end{lemma}

\proof
For all absolutely continuous functions $f$ for which the expectations below exist, 
\beas
\sigma^2 \E f'(Y^*) &=& \sigma^2 \E f'(UY^\dagger + (1-U)Y^\ddagger) \\
           &=& \sigma^2 \E \left(\frac{f(Y^\dagger)-f(Y^\ddagger)}{Y^\dagger-Y^\ddagger}\right) \\
					 &=& \frac{\sigma^2}{2\left(\lambda \sigma^2-\E Y'R \right)}\E \left(\left(\frac{f(Y'')-f(Y')}{Y''-Y'}\right)(Y''-Y')^2\right)\\
					 &=& \frac{\sigma^2}{2\left(\lambda \sigma^2 -\E Y'R \right)} \E((f(Y'')-f(Y'))(Y''-Y')) \\
					 &=& \frac{\sigma^2}{\lambda \sigma^2 - \E Y'R} \E(Y'f(Y')-Y''f(Y')) \\
					 &=& \frac{\lambda \sigma^2}{\lambda \sigma^2 - \E Y'R} \E(Y'f(Y')) - \frac{\sigma^2}{\lambda \sigma^2-\E Y'R} \E Rf(Y'),
\enas
where we have used \eqref{diff2} and \eqref{approxsteinpair2} in the third and the last equalities, respectively.

Thus
\beas
\E Y'f(Y') &=& \frac{\lambda \sigma^2 - \E Y'R}{\sigma^2 \lambda} \sigma^2 \E f'(Y^*) + \frac{\E Rf(Y')}{\lambda} \\
         &=& \sigma^2 \E f'(Y^*) - \frac{\E Y'R}{\lambda}\E f'(Y^*) + \frac{\E Rf(Y')}{\lambda}.
\enas
\bbox

\begin{remark}
One may notice that we construct an approximate zero bias coupling through an exchangeable pair. An important reason that we develop this coupling technique instead of simply using Stein's method of exchangeable pairs is that we aim to avoid the calculation of the term $\Var\left(\E[(Y''-Y')^2|Y']\right)$ that can be difficult to compute in many cases. The reader will see an example in Section \ref{comb} that we take one more step that might not be very easy to construct an approximate zero bias coupling but all the computations after that are straightforward. 
\end{remark}

Next the following result shows how to construct an approximate zero bias distribution of $Y'$ using an approximate Stein pair $(Y',Y'')$, when the latter is a function of some underlying random variables $\xi_{\alpha}, \alpha \in \chi$ and a random index $\mathbf{I}$. It is a minor variation of Lemma 4.4 of \cite{CGS11}, with a $\lambda$-Stein pair and $2\lambda \sigma^2$ there respectively replaced by an approximate $\lambda,R$-Stein pair and $\E(Y'-Y'')^2=2\left(\lambda \sigma^2 -\E Y'R \right)$ here. The proof is omitted, being similar under these replacements.

\begin{lemma} \label{Wdagger2}
Let $F(y',y'')$ be the distribution of an approximate $\lambda,R$-Stein pair $(Y',Y'')$ and suppose there exist a distribution
\bea \label{distixi}
F({\bf i},\xi_{\alpha},\alpha \in \chi)
\ena
and an $\mathbb{R}^2$ valued function $(y',y'') = \psi({\bf i},\xi_{\alpha},\alpha \in \chi)$ such that when ${\bf I}$ and $\{\Xi_{\alpha}, \alpha \in \chi \}$ have distribution \eqref{distixi} then
\beas
(Y',Y'') = \psi({\bf I},\Xi_{\alpha}, \alpha \in \chi )
\enas
has distribution $F(y',y'')$. If ${\bf I}^{\dagger}$, $\{\Xi^\dagger_{\alpha}, \alpha \in \chi\}$ have distribution
\bea \label{squarebias}
dF^\dagger({\bf i},\xi_{\alpha}, \alpha \in \chi) = \frac{(y'-y'')^2}{\E (Y'-Y'')^2} dF({\bf i},\xi_{\alpha}, \alpha \in \chi)
\ena
then the pair 
\beas
(Y^\dagger,Y^\ddagger) = \psi({\bf I}^{\dagger},\Xi^\dagger_{\alpha}, \alpha \in \chi)
\enas
has distribution $F^\dagger(y^\dagger,y^\ddagger)$ satisfying
\bea \label{eq:square.bias.pair}
dF^\dagger(y',y'') = \frac{(y'-y'')^2}{\E (Y'-Y'')^2}dF(y',y'').
\ena
\end{lemma}

In the following, for functions $f:\mathbb{R}\rightarrow \mathbb{R}$, we let $|f|_\infty = \sup_{x\in \mathbb{R}}|f(x)|$ be the supremum norm. Theorems \ref{thm:approx.zb.dist.l1} and \ref{thm:approx.zb.dist.Kol} below provide respectively $L^1$ and $L^\infty$ bounds between $Y'$ and a standard normal random variable in term of $Y^*$ satisfying \ref{approxzerobias}.

\begin{theorem} \label{thm:approx.zb.dist.l1}
Let $Y'$ be a mean zero, variance $\sigma^2>0$ random variable, and $Y^*$ and $R(Y')$ be defined on the same space as $Y'$, satisfying \eqref{approxzerobias}, then
\bea \label{l1approx}
d_1(\mathcal{L}(Y'/\sigma),\mathcal{L}(Z)) \le \frac{2}{\sigma}\E|Y^*-Y'| + \sqrt{\frac{2}{\pi}}\frac{|\E Y'R|}{\sigma^2 \lambda} + \frac{2\E |R|}{\sigma \lambda},
\ena
where $Z$ is a standard normal random variable.
\end{theorem}

\proof
For given $h \in \mathbb{L}$ let $f$ be the unique bounded solution to the Stein equation
\beas
f'(w)-w f(w) = h(w)-Nh.
\enas
Then, (see e.g. \cite{CGS11} Lemma 2.4),
\beas
|f|_{\infty} \le 2,|f'|_{\infty} \le \sqrt{\frac{2}{\pi}} \qmq{and} |f''|_{\infty} \le 2.
\enas
Letting $g(x) = f(x/\sigma)$, we have $g'(x) = (1/\sigma)f'(x/\sigma)$ , and thus
\beas
&&|\E h(Y'/\sigma) - Nh| = \left|\E \left(f'\left(\frac{Y'}{\sigma}\right)-\frac{Y'}{\sigma}f\left(\frac{Y'}{\sigma} \right) \right) \right| \nn\\
                       &&\hspace{10pt}= \left|\E \left(f'\left(\frac{Y'}{\sigma}\right)-\frac{1}{\sigma}Y'g\left(Y' \right) \right) \right| \\
&&\hspace{10pt} =\left|\E \left(f'\left(\frac{Y'}{\sigma} \right)-\sigma g'(Y^*)\right)+ \frac{\E Y'R}{\sigma\lambda}\E g'(Y^*) - \frac{\E Rg(Y')}{\sigma \lambda}\right| \\
&&\hspace{10pt}\le \E \left|f'\left(\frac{Y'}{\sigma} \right) - f'\left(\frac{Y^*}{\sigma} \right) \right|+  \left| \frac{\E Y'R}{\sigma^2\lambda}\E f'\left(\frac{Y^*}{\sigma} \right)\right| +  \left|\frac{\E Rf(Y'/\sigma)}{\sigma\lambda}\right| \\
&&\hspace{10pt}\le \frac{2}{\sigma}\E |Y^*-Y'| + \sqrt{\frac{2}{\pi}} \frac{|\E Y'R|}{\sigma^2\lambda} + \frac{2\E |R|}{\sigma\lambda},
\enas
where we have applied \eqref{approxzerobias} with $f$ replaced by $g$ in the third equality.

\bbox

\begin{theorem} \label{thm:approx.zb.dist.Kol}
Let $Y'$ be a mean zero, variance $\sigma^2>0$ random variable, and $Y^*$ and $R(Y')$ be defined on the same space as $Y'$, satisfying \eqref{approxzerobias} and $|Y^*-Y|\le\delta$. Then
\bea \label{Kolapprox}
d_\infty(\mathcal{L}(Y'/\sigma),\mathcal{L}(Z)) \le \frac{\delta(1+1/\sqrt{2\pi}+\sqrt{2\pi}/4)}{\sigma} + \frac{|\E Y'R|}{\sigma^2\lambda} + \frac{\sqrt{2\pi}\E|R|}{4\sigma\lambda},
\ena
where $Z$ is a standard normal random variable.
\end{theorem}

\proof We follow the proof of Theorem 5.1 of \cite{CGS11} which obtained the the same type of bound using zero biasing.
Let $z \in \mathbb{R}$, $\epsilon = \delta/\sigma$ and $f$ be the solution of the equation
\beas
f'(w)-wf(w) = \mathbf{1}_{ \{w \le z-\epsilon \}} - P(Z \le z-\epsilon).
\enas
Then, (see e.g. \cite{CGS11} Lemma 2.3),
\bea \label{steineq:bounds2}
|f|_{\infty} \le \frac{\sqrt{2\pi}}{4},|f'|_{\infty} \le 1,
\ena
and for all $w$, $u$ and $v$,
\bea \label{steineq:bounds3}
|(w+u)f(w+u) - (w+v)f(w+v)| \le (|w|+\sqrt{2\pi}/4)(|u|+|v|).
\ena

Then we have
\beas
f'(Y^*/\sigma) - (Y^*/\sigma)f(Y^*/\sigma) &=& \mathbf{1}_{ \{Y^*/\sigma \le z-\epsilon \}} - P(Z \le z-\epsilon) \\  
                                           &\le& \mathbf{1}_{ \{Y'/\sigma \le z \}} - P(Z \le z-\epsilon).
\enas
Letting $g(x) = f(x/\sigma)$, we have $g'(x) = (1/\sigma)f'(x/\sigma)$, and
\beas
P(Y'/\sigma \le z) - P\left(Z < z\right) &=& (P(Z < z-\epsilon)-P(Z < z)) \\ 
                                             && \hspace{10pt}+P(Y'/\sigma \le z)-P(Z < z-\epsilon) \\
                    &\ge& -\epsilon/\sqrt{2\pi} + P(Y'/\sigma \le z)-P(Z < z-\epsilon) \\
										&\ge& -\epsilon/\sqrt{2\pi}  + \E[f'(Y^*/\sigma)-(Y^*/\sigma)f(Y^*/\sigma)] \\
										&=& -\epsilon/\sqrt{2\pi} + (1/\sigma)\E[\sigma^2 g'(Y^*)-Y^* g(Y^*)].
\enas
Then, using \eqref{approxzerobias} and applying \eqref{steineq:bounds2} in the last inequality,
\bea \label{kol}
P(Y'/\sigma \le z) - P\left(Z < z\right) &\ge& -\frac{\epsilon}{\sqrt{2\pi}} + \frac{1}{\sigma}\E\left[Y'g(Y')-Y^*g(Y^*)\right] \nn \\
&&\hspace{10pt} + \frac{\E Y'R}{\sigma\lambda}\E g'(Y^*) - \frac{\E Rg(Y')}{\sigma\lambda} \nn\\
&=& -\frac{\epsilon}{\sqrt{2\pi}} + \E\left[\frac{Y'}{\sigma}f\left(\frac{Y'}{\sigma}\right)-\frac{Y^*}{\sigma}f\left(\frac{Y^*}{\sigma}\right)\right] \nn\\
           && \hspace{10pt}+ \frac{\E Y'R}{\sigma^2\lambda}\E f'\left(\frac{Y^*}{\sigma}\right) - \frac{\E Rf(Y'/\sigma)}{\sigma\lambda}  \nn\\
&\ge&  -\frac{\epsilon}{\sqrt{2\pi}} + \E\left[\frac{Y'}{\sigma}f\left(\frac{Y'}{\sigma}\right)-\frac{Y^*}{\sigma}f\left(\frac{Y^*}{\sigma}\right)\right] \nn \\
        && \hspace{10pt} - \frac{|\E Y'R|}{\sigma^2\lambda} - \frac{\sqrt{2\pi}\E |R|}{4\sigma\lambda} .  
\ena
Writing $\Delta = Y^*/\sigma - Y'/\sigma$ and applying \eqref{steineq:bounds3} and that $|\Delta| \le \epsilon$ yields
\beas
\left|\E\left[\frac{Y'}{\sigma}f\left(\frac{Y'}{\sigma}\right)-\frac{Y^*}{\sigma}f\left(\frac{Y^*}{\sigma}\right)\right]\right| 
     &=& \left|\E\left[\frac{Y'}{\sigma}f\left(\frac{Y'}{\sigma}\right)-\left(\frac{Y'}{\sigma}+\Delta\right)f\left(\frac{Y'}{\sigma}+\Delta\right)\right]\right|  \\
		 &\le& \E \left[\left(\frac{|Y'|}{\sigma}+\frac{\sqrt{2\pi}}{4}\right)|\Delta|\right] \\
		 &\le& \epsilon(1+\sqrt{2\pi}/4) = \frac{\delta}{\sigma}(1+\sqrt{2\pi}/4).
\enas
Using this inequality in \eqref{kol} yields
\beas
P(Y'/\sigma \le z) - P\left(Z < z\right) \ge -\frac{\delta}{\sigma}\left(\frac{1}{\sqrt{2\pi}}+1+\frac{\sqrt{2\pi}}{4}\right) - \frac{|\E Y'R|}{\sigma^2\lambda} - \frac{\sqrt{2\pi}\E |R|}{4\sigma\lambda}.
\enas
A similar argument yields the reverse inequality.

\bbox

We end up this section by mentioning a connection between our work in this section and the known result when $R=0$. 
\begin{remark}
When $R=0$, the construction of Lemma \ref{Wdagger} yields the construction of the zero bias distribution as in Lemma \ref{zerobiaslemma}. Furthermore, the $L^1$ and $L^\infty$ bounds in Theorems \ref{thm:approx.zb.dist.l1} and \ref{thm:approx.zb.dist.Kol} reduce to the bounds in \eqref{oldl1} and \eqref{oldKol}, respectively.
\end{remark}

\section{Ewens measure} \label{Ewens}
In this section, we briefly describe the Ewens measure and state some necessary properties. Let $\mathcal{S}_n$ denotes the symmetric group. The Ewens distribution ${\cal E}_\theta$ on the symmetric group ${\cal S}_n$ with parameter $\theta > 0$, was first introduced in \cite{Ewe72} and used in population genetics to describe the probabilities associated with the number of times that different alleles are observed in the sample; see also \cite{ABT03} for the description in mathematical context. In the following, we let $\mathbb{N}_k=[k,\infty) \cap \mathbb{Z}$ and for $x \in \mathbb{R}$, $n \in \mathbb{N}_1$, we use the notations
\beas
x^{(n)} = x(x+1)\cdots(x+n-1)  \text{ \ \ and \ \ } x_{(n)} = x(x-1)\cdots(x-n+1).
\enas
Given a permutation $\pi \in \mathcal{S}_n$, the Ewens measure is given by 
\bea \label{ewensdef}
P_{\theta}(\pi) = \frac{\theta^{\#(\pi)}}{\theta^{(n)}},
\ena
where $\#(\pi)$ denotes the number of cycles of $\pi$. We note that ${\cal E}_\theta$ specializes to the uniform distribution over all permutations when $\theta=1$.

The Ewens measure ${\cal E}_\theta$ can be defined equivalently in term of $c_1(\pi),\ldots,c_n(\pi)$ as follows,
\bea \label{ewensdef2}
P_\theta(c_1,\ldots,c_n) =\mathbf{1}\left\{\sum_{j=1}^n jc_j = n\right\}\frac{n!}{\theta^{(n)}}\prod_{j=1}^n\frac{\theta^{c_j}}{j^{c_j}c_j!},
\ena
where $c_q(\pi)$ is the number of $q$ cycles of $\pi$ and we write $c_q$ for $c_q(\pi)$ for simplicity.

A permutation $\pi_n \in \mathcal{S}_n$ with the distribution ${\cal E}_\theta$ can be constructed by the `so called' the Chinese restaurant process (see e.g. \cite{Ald85} and \cite{Pit96}), as follows. For $n=1$, $\pi_1$ is the unique permutation that maps $1$ to $1$ in $\mathcal{S}_1$. For $n \ge 2$, we construct $\pi_n$ from $\pi_{n-1}$ by either adding $n$ as a fixed point with probability $\theta/(\theta+n-1)$, or by inserting $n$ uniformly into one of $n-1$ locations inside a cycle of $\pi_{n-1}$, so each with probability $1/(\theta+n-1)$.

For $\sigma$ a permutation of the elements of $[n]$ and $B \subset [n]$, we will consider the reduced permutation $\sigma\setminus B$ of the elements $[n]\setminus B$ whose cycle representation is obtained by deleting all elements of $B$ in the cycle representation of $\sigma$. For instance, if $n=5$ and the cycle representation of $\sigma$ is $(1)(2435)$ and $B=\{1,2\}$ then $\sigma\setminus B$ has representation $(354)$. Also, let $\sigma_B$ be the permutation whose cycle structure is obtained by taking the cycle structure of $\sigma$ and removing all cycles that contain any element of $B^c$. Here, for instance,  $\sigma_B$ has cycle structure $(1)$. With $\# (\tau)$ denoting the number of cycles of the permutation $\tau$, we easily see that $\#(\sigma)=\#(\sigma/B) +\#(\sigma_B)$, as any cycle of $\sigma$ either contains, or does not contain, some element of $B^c$.

For $B \subset [n]$, Propositions \ref{probEwens} and \ref{conditionEwens}
	that follow respectively provide the joint unconditional probability that $\pi(i) = \xi_i,i \in B$, and conditional probability that $\pi(i) = \xi_i,i \in B^c$ given $\{\pi(i),i \in B\}$ under the Ewens distribution.

\begin{proposition} \label{probEwens}
Let $\pi$ be a permutation of $[n]$ with distribution $\mathcal{E}_\theta$, $\theta>0$ and $B \subseteq [n]$. Then, for $\xi_i,i\in B$ distinct elements of $[n]$,
\beas
P(\pi(i)=\xi_i,i \in B) = \frac{\theta^{\#(\pi_B)}}{(\theta+n-1)_{(|B|)}}.
\enas 
\end{proposition}
\proof 
We prove this lemma by induction on the size of $B$. Since it is clear by \eqref{ewensdef} that the distribution of $\pi$ depends only on the number of cycles, it is sufficient to prove the result for $B_m = \{n-m+1,n-m+2,\ldots,n \}$ and $m \in [n]$. For $B_1 = \{n \}$, by the starting configuration in the construction of $\pi$ via the Chinese restaurant process described above, we immediately have 
\beas
P_\theta(\pi(n) = k) =
\begin{cases}
\frac{1}{\theta+n-1} \text{ \ if \ } k \ne n,\\
\frac{\theta}{\theta+n-1} \text{ \ if \ } k = n
\end{cases}
= \frac{\theta^{\#(\pi_{B_1})}}{(\theta+n-1)_{(|B_1|)}}.
\enas
Now we assume that the claim is true for $B_{m-1}$. To prove the result for $B_m$, we recall that the Chinese restaurant process either adds $n-m+1$ as a fixed point with probability $\theta/(\theta+n-m)$, or inserts $n-m+1$ uniformly into one of $n-m$ locations inside a cycle of $\pi_{n-m}$, so each with probability $1/(\theta+n-m)$. Hence, using the assumption that the result holds for $B_{m-1}$, we have
\beas
P_\theta(\pi(i) = \xi_i,i \in B_{m}) &=&
\begin{cases}
\frac{1}{\theta+n-m}\frac{\theta^{\#(\pi_{B_{m-1}})}}{(\theta+n-1)_{(m-1)}}  \text{ \ if $n-m+1$ is not a fixed point  } ,\\
\frac{\theta}{\theta+n-m}\frac{\theta^{\#(\pi_{B_{m-1}})}}{(\theta+n-1)_{(m-1)}} \text{ \ if $n-m+1$ is a fixed point } 
\end{cases} \\
&=& \frac{\theta^{\#(\pi_{B_m})}}{(\theta+n-1)_{(m)}}=\frac{\theta^{\#(\pi_{B_m})}}{(\theta+n-1)_{(|B_m|)}}.
\enas
\bbox



\begin{proposition}  \label{conditionEwens}
Let $\pi$ be a permutation of $[n]$ with distribution $\mathcal{E}_\theta$, $\theta>0$ and $B \subset [n]$. Then, for $\xi_i,i\in [n]$ distinct elements of $[n]$,
\beas
P(\pi(i)=\xi_i,i \in B^c| \pi(k) = \xi_k,k\in B) = \frac{\theta^{\#(\pi  \setminus B)}}{\theta^{(|B^c|)}}.
\enas 
\end{proposition} 
\proof
Using the definition of conditional probability and \eqref{ewensdef} and applying Proposition \ref{probEwens}, we have
\beas
P(\pi(i)=\xi_i,i \in B^c| \pi(k) = \xi_k,k\in B) &=& \frac{P(\pi(i)=\xi_i,i \in [n])}{P(\pi(k) = \xi_k,k\in B)} \\
                                                 &=& \frac{\theta^{\#(\pi )}/\theta^{(n)}}{\theta^{\#(\pi_B)}/(\theta+n-1)_{(|B|)}} \\
																								 &=& \frac{\theta^{\#(\pi   \setminus B)}}{\theta^{(|B^c|)}},
\enas
where we have used $\#(\pi \setminus B) = \#(\pi )-\#(\pi_B) $ in the last equality.

\bbox

The joint moments of $c=(c_1,\ldots,c_n)$ in the uniform case were established in \cite{Wat74} (See also \cite{ABT03}). Using the similar argument, in the following proposition, we generalize the result to the joint moments of $c=(c_1,\ldots,c_n)$ under the distribution $\mathcal{E}_\theta$ for any $\theta>0$. Note that, with $\theta=1$, the proposition below is exactly the same as the result in \cite{Wat74} and \cite{ABT03}.  

\begin{proposition} \label{c1moment}
Let $c=(c_1,\ldots,c_n)$ be a cycle type of $\pi \in \mathcal{S}_n$ with distribution $\mathcal{E}_\theta$ with $\theta > 0$. Then, for $m_1,\ldots,m_n \in \mathbb{N}_0$ with $m = \sum_{j=1}^n jm_j$,
\beas
\E\left(\prod_{j=1}^n (c_j)_{(m_j)}\right) = \left(\frac{n_{(m)}}{(\theta+n-1)_{(m)}} \prod_{j=1}^n \left(\frac{\theta}{j}\right)^{m_j} \right)
                                                              \mathbf{1}\left\{\sum_{j=1}^n jc_j \le n\right\}.																												
\enas
\end{proposition}
\proof
Using $(c_j)_{(m_j)}/c_j! = 1/(c_j-m_j)!$ when $c_j \ge m_j$, we have
\beas
&&\E\left(\prod_{j=1}^n (c_j)_{(m_j)}\right) = \sum_{c \in \mathbb{N}_0^n} P_\theta(c) \prod_{j=1}^n (c_j)_{(m_j)} \\
  && \hspace{10pt}= \sum_{c:c_j\ge m_j} \mathbf{1}\left\{\sum_{j=1}^n jc_j = n\right\} \frac{n!}{\theta^{(n)}}\prod_{j=1}^n \frac{(c_j)_{(m_j)} \theta^{c_j}}{j^{c_j}c_j!} \\
	&&\hspace{10pt}= \frac{n_{(m)}}{(\theta+n-1)_{(m)}} \prod_{j=1}^n \left(\frac{\theta}{j}\right)^{m_j}  
	     \sum_{d \in \mathbb{N}_0^n} \mathbf{1}\left\{\sum_{j=1}^n jd_j = n-m\right\} \frac{(n-m)!}{\theta^{(n-m)}}\prod_{j=1}^n \frac{\theta^{d_j}}{j^{d_j}d_j!} \\
  &&\hspace{10pt}= \frac{n_{(m)}}{(\theta+n-1)_{(m)}} \prod_{j=1}^n \left(\frac{\theta}{j}\right)^{m_j}  
	   \sum_{d \in \mathbb{N}_0^n} \mathbf{1}\left\{\sum_{j=1}^{n-m} jd_j = n-m\right\} \frac{(n-m)!}{\theta^{(n-m)}}\prod_{j=1}^{n-m} \frac{\theta^{d_j}}{j^{d_j}d_j!}\\ 
	&&\hspace{10pt}= \left(\frac{n_{(m)}}{(\theta+n-1)_{(m)}} \prod_{j=1}^n \left(\frac{\theta}{j}\right)^{m_j}\right) \mathbf{1}\left\{\sum_{j=1}^n jc_j \le n\right\}
\enas
where $d_j$ corresponds to $c_j-m_j$ and the sum over $d$ in the second last line is one as it is taken over all possibilities of $\mathcal{S}_{n-m}$.

\bbox


\section{Combinatorial CLT under the Ewens measure} \label{comb}
In this section, following Section 6.1 of \cite{CGS11}, we study the distribution, introduced in \cite{Hoe51}, of 
\bea \label{combdef}
Y =\sum_{i=1}^n a_{i,\pi(i)}
\ena
where $A \in \mathbb{R}^{n \times n}$ is a given real matrix with components $\{a_{i,j} \}_{i,j=1}^n$ and $\pi \in \mathcal{S}_n$ has the Ewens distribution.

A distribution on $\mathcal{S}_n$ is said to be constant on cycle type if the probability of any permutation $\pi \in {\cal S}_n$ depends only on the cycle type $(c_1,\ldots,c_n)$. The work \cite{CGS11} studied the distribution of
\eqref{combdef} where $\pi \in \mathcal{S}_n$ has distribution constant on cycle type with no fixed points. It follows from \eqref{ewensdef2} directly that the Ewens distribution is constant on cycle type and allows fixed point. Therefore, though several techniques in Section 6.1.2 of \cite{CGS11} apply here, the main proofs and the coupling construction do not.

Letting
\bea \label{adotdot}
a_{\bullet,\bullet} = \frac{1}{n(\theta+n-1)} \left(\theta \sum_{i=1}^n a_{i,i}+ \sum_{i \ne j}a_{ij}\right),
\ena 
applying Proposition \ref{probEwens} with $B = \{ i\}$, we have
\bea \label{meanY}
\E Y = \sum_{i=1}^n \E a_{i,\pi(i)} =  \sum_{i=1}^n \left(\frac{\theta}{\theta+n-1}a_{i,i}+\sum_{j,j\ne i} \frac{1}{\theta+n-1}a_{i,j}\right) = n a_{\bullet,\bullet}.
\ena
Letting
\beas
\widehat{a}_{i,j} = a_{i,j} - a_{\bullet,\bullet},
\enas
and using \eqref{meanY}, we have 
\beas
\E \left[\sum_{i=1}^n \widehat{a}_{i,\pi(i)}\right]=\E \left[\sum_{i=1}^n (a_{i,\pi(i)}-a_{\bullet,\bullet})\right] = \E \left[\sum_{i=1}^n a_{i,\pi(i)} - na_{\bullet,\bullet}\right] = 0.
\enas
As a consequence, replacing $a_{i,j}$ by $\widehat{a}_{i,j}$, we may without loss of generality assume that
\bea \label{zeroassumption}
\E Y = 0 \text{ \ and \ } a_{\bullet,\bullet} = 0, 
\ena 
and for simplicity, as in \cite{CGS11}, we consider only the symmetric case, that is, $a_{i,j} = a_{j,i}$ for all $i,j \in [n]$.

To rule out trivial cases, we assume in what follows that $\sigma^2>0$. We later calculate $\sigma^2$ explicitly in \eqref{variance} of Lemma \ref{combvar} and discuss in Remark \ref{sigmaorder} that it is of order $n$ when the elements of $A$ are well chosen in some sense. In this case, there exists $N \in \mathbb{N}_1$ such that $\sigma^2>0$ for $n > N$.

In the following theorem, we obtain upper bounds for the $L^1$ and $L^\infty$ distances between $Y$ given in \eqref{combdef} with the distribution $\mathcal{E}_\theta$ and a standard normal random variable and lower bounds for the $L^\infty$ distance in the special case that the matrix $A$ is integer-valued. Below, we consider the symbols $\pi$ and $Y$ interchageable with $\pi'$ and $Y'$, respectively.

\begin{theorem} \label{L1theorem}
Let $n \ge 6$ and $\{a_{i,j} \}_{i,j=1}^n$ be an array of real numbers satisfying
\beas
a_{i,j} = a_{j,i}.
\enas
Let $\pi' \in \mathcal{S}_n$ be a random permutation with the distribution ${\cal E}_\theta$, with $\theta > 0$.
Then, with $Y'$ the sum in \eqref{combdef} with $\pi'$ replacing $\pi$, assuming $\Var(Y') = \sigma^2>0$, and letting $W = (Y' -\E Y')/\sigma$ and $Z$ a standard normal random variable, 
\bea \label{mainbound}
d_1(\mathcal{L}(W),\mathcal{L}(Z)) \le  \frac{\alpha_1(\theta,M,n)}{\sigma},
\ena
and
\bea \label{mainbound2}
d_\infty(\mathcal{L}(W),\mathcal{L}(Z)) \le  \frac{\alpha_2(\theta,M,n)}{\sigma},
\ena
where 
\beas
\alpha_1(\theta,M,n) &=&  40M + \kappa_{\theta,n,1}\sqrt{2/\pi}M\left(3 + \frac{\theta+1}{n-1}\right)+\frac{\kappa_{\theta,n,2}\sqrt{2/\pi}M}{n-1} \nn\\
                       &&  + \theta M \left(1.2\sqrt{2/\pi}+\frac{1.2(6n+4\theta-5)}{\theta+n-1} + \frac{\theta n}{(\theta+n-1)_{(2)}}\right),
\enas
and
\beas
\alpha_2(\theta,M,n) &=&  20(1+1/\sqrt{2\pi}+\sqrt{2\pi}/4) M  + \kappa_{\theta,n,1}M\left(3 + \frac{\theta+1}{n-1}\right)+\frac{\kappa_{\theta,n,2}M}{n-1} \nn\\
    && + \theta M \left(1.2+\frac{0.15\sqrt{2\pi}(6n+4\theta-5)}{\theta+n-1} + \frac{0.125\sqrt{2\pi}\theta n}{(\theta+n-1)_{(2)}}\right)
\enas
with 
\beas
M = \max_{i,j} |a_{i,j} -a_{\bullet,\bullet}|, \,  a_{\bullet,\bullet} \text{ \ as in \eqref{adotdot}},
\enas
\bea \label{kappadef}
\kappa_{\theta,n,1} = \sqrt{\frac{\theta^2n_{(2)}}{(\theta+n-1)_{(2)}} + \frac{\theta n}{\theta+n-1}}
\ena
and
\bea \label{kappadef2}
\kappa_{\theta,n,2} = \sqrt{\frac{\theta^4n_{(4)}}{(\theta+n-1)_{(4)}}+\frac{4\theta^3 n_{(3)}}{(\theta+n-1)_{(3)}}+\frac{2\theta^2n_{(2)}}{(\theta+n-1)_{(2)}}}.
\ena
In particular, if $a_{i,j}$, $i,j \in [n]$ are all integers, then 
\bea \label{mainlowerbound}
\frac{1}{6\sqrt{3} \sigma+3} \le d_\infty(\mathcal{L}(W),\mathcal{L}(Z)) \le  \frac{\alpha_2(\theta,M,n)}{\sigma}.
\ena
\end{theorem}

Remark \ref{thetaremark} below discusses the behavior of \eqref{kappadef} and \eqref{kappadef2} in $\theta$ and the bounds on the rates of convergence in $n$.

\begin{remark} \label{thetaremark}
\begin{enumerate} 
  \item $\kappa_{\theta,n,1}$ and $\kappa_{\theta,n,2}$  given in \eqref{kappadef} and \eqref{kappadef2} are of constant order in $n$ since
	\beas
	\kappa_{\theta,n,1} \rightarrow \sqrt{\theta^2+\theta} \text{ \ \ and \ \ } \kappa_{\theta,n,2} \rightarrow \sqrt{\theta^4+4\theta^3+2\theta^2} \text{ \ \ as \ \ } n \rightarrow \infty. 
	\enas
	In particular, $\kappa_{1,n,1} = \sqrt{2}$ and $\kappa_{1,n,2} = \sqrt{7}$.
	\item Since $\kappa_{\theta,n,1}$ and $\kappa_{\theta,n,2}$ are of constant order in $n$, $\alpha_1(\theta,M,n)$ and $\alpha_2(\theta,M,n)$ are also of constant order in $n$ for all $\theta > 0$. Therefore, the $L^1$ and $L^\infty$ bounds on the right hand side of \eqref{mainbound} and \eqref{mainbound2} are of order $\sigma^{-1}$ which is of the same order as the $L^1$ bound in the uniform case in Theorem 4.8 of \cite{CGS11}. The uniform distribution corresponds to the special case of the Ewens distribution with $\theta = 1$ and the result in this case is presented in Corollary \ref{uniform} below. The order of $\sigma^2$ in $n$ is considered in Remark \ref{sigmaorder}, where we find that if $A=[a_{i,j}]$ is chosen well in some sense then $\sigma^2$ will be of order $n$, implying that the bounds in \eqref{mainbound} and \eqref{mainbound2} are of order $n^{-1/2}$. 
	\item In the case that $a_{i,j}$, $i,j \in [n]$ are interger-valued, the $L^\infty$ upper and lower bounds in \eqref{mainlowerbound} are of the same order $\sigma^{-1}$,  which is thus the optimal order for the $L^\infty$ distance. 
	\item It is easy to see that $\alpha_1(\theta,M,n)$ and $\alpha_2(\theta,M,n)$ are increasing in $\theta$ and thus so are the $L^1$ and $L^\infty$ bounds in \eqref{mainbound} and \eqref{mainbound2}. In fact, the expressions $\alpha_1(\theta,M,n)$ and $\alpha_2(\theta,M,n)$ are obtained from \eqref{l1approx} and \eqref{Kolapprox}, respectively, which depend on  $R$. This remainder $R$, given explicitly in Lemma \ref{combapproxstein} below, depends on the number of fixed points of $\pi$, which become more likely as $\theta$ increases. 
	\item As $\theta$ tends to zero,
\beas
\alpha_1(\theta,M,n)\rightarrow 40M \text{ \ \ and \ \ } \alpha_2(\theta,M,n) \rightarrow 20(1+1/\sqrt{2\pi}+\sqrt{2\pi}/4) M.
\enas 
Therefore, the $L^1$ and $L^\infty$ bounds on the right hand side of \eqref{mainbound} and \eqref{mainbound2} converge to $40M/\sigma$ and  $20(1+1/\sqrt{2\pi}+\sqrt{2\pi}/4) M/\sigma$, respectively. This corresponds to the case where a large number of cycles is unlikely.
\end{enumerate}

\end{remark} 

Next we present Corollary \ref{uniform} that specializes Theorem \ref{L1theorem} to the case where the random permutation $\pi'$ has the uniform distribution, corresponding to the special case of the Ewens distribution with $\theta = 1$. Indeed, the result immediately follows from Theorem \ref{L1theorem} by applying the bounds
\beas
\alpha_1(1,M,n) =  \left(47.2+6/\sqrt{\pi}+1.2\sqrt{2/\pi}+\frac{5+\sqrt{14}}{\sqrt{\pi}(n-1)}-\frac{1.2}{n}\right)M  \le 53M         
\enas
and 
\beas
\alpha_2(1,M,n) &=& \Bigg(21.2+3\sqrt{2}+5.9\sqrt{2\pi}+20/\sqrt{2\pi} \\
                 && \hspace{50pt}+ \frac{2\sqrt{2}+\sqrt{7}+0.125\sqrt{2\pi}}{n-1}-\frac{0.15\sqrt{2\pi}}{n}\Bigg)M 
                          \le 50M,
\enas
which hold for all $n \ge 6$.

\begin{corollary} \label{uniform}
Let $n \ge 6$ and $\{a_{i,j} \}_{i,j=1}^n$ be an array of real numbers satisfying
\beas
a_{i,j} = a_{j,i}.
\enas
Let $\pi'$ be a random permutation with the uniform distribution over ${\cal S}_n$.
Then, with $Y'$ the sum in \eqref{combdef} with $\pi'$ replacing $\pi$, assuming $\Var(Y') = \sigma^2>0$, and letting $W = (Y' -\E Y')/\sigma$ and $Z$ a standard normal random variable,
\beas
d_1(\mathcal{L}(W),\mathcal{L}(Z)) \le  \frac{53M }{\sigma} ,
\enas
and
\beas
d_\infty(\mathcal{L}(W),\mathcal{L}(Z)) \le \frac{50M}{\sigma},
\enas
where 
\beas
M = \max_{i,j} |a_{i,j} -a_{\bullet,\bullet}| \text{ \ and \ } \  a_{\bullet,\bullet} \text{ \ is as in \eqref{adotdot}}.
\enas
In particular, if $a_{i,j}$, $i,j \in [n]$ are all integers, then
\beas
\frac{1}{6\sqrt{3} \sigma+3} \le d_\infty(\mathcal{L}(W),\mathcal{L}(Z)) \le  \frac{50M}{\sigma}.
\enas
\end{corollary}

Section 6.1.2 of \cite{CGS11} proved the main $L^1$ and $L^\infty$ bounds by first considering each cycle type $c:=(c_1,\ldots,c_n)$ separately and then combining all cases. Since fixed points are allowed in the present work, $P_{\theta}(\pi(i)=i|c)$ varies as $c$ changes. This dependence on $c$ makes it difficult to follow the same proof structure and coupling construction as there. As a result, we construct a new coupling for the proof of Theorem \ref{L1theorem}. 

To prove the main results of this section, our target is therefore to apply Theorems \ref{thm:approx.zb.dist.l1} and \ref{thm:approx.zb.dist.Kol}. Hence we will first construct an approximate $\lambda,R$-Stein pair and then an approximate zero bias coupling through the helps of Lemmas \ref{Wdagger} and \ref{Wdagger2}. For this purpose, we present a sequence of lemmas below. The first two lemmas were proved in \cite{CGS11}. The proofs of Lemmas \ref{combapproxstein} and \ref{combvar}, though important, the first closely follows the proof of Lemma 6.9 of \cite{CGS11} and the latter is straightforward, can be found in the Appendix. 
Lemma \ref{partition} forms a partition of the space based on the cycle structure of $\pi \in \mathcal{S}_n$. Using that partition, Lemma \ref{diff} expresses the difference in the values taken on by the exchangeable pair coupling given in Lemma \ref{combapproxstein}. Below, for $i,j \in [n]$, we write $i\sim j$ if $i$ and $j$ are in the same cycle, let $|i|$ be the length of the cycle containing $i$ and let $\tau_{i,j}$, $i,j \in [n]$ be the permutation that transposes $i$ and $j$.

\begin{lemma}[\cite{CGS11}] \label{partition} 
Let $\pi$ be a fixed permutation. For any $i \ne j$, distinct elements of $[n]$, the sets $A_0,\ldots,A_5$ form a partition of the space where,
\beas
&& A_0 = \{|\{i,j,\pi(i),\pi(j) \}|=2 \} \\
&& A_1 = \{|i|=1, |j| \ge 2 \}, \ \ \ A_2 = \{|i| \ge 2, |j|=1 \} \\
&& A_3 = \{|i| \ge 3, \pi(i) = j \}, \ \ \ A_4 = \{|j| \ge 3, \pi(j) = i \} \text{ \ \ and } \\
&& A_5 = \{|\{i,j,\pi(i),\pi(j) \}|=4 \}.
\enas
Additionally, the sets $A_{0,1}$ and $A_{0,2}$ partition $A_0$ where
\beas
A_{0,1} = \{\pi(i)=i,\pi(j)=j \}, \ \ \ A_{0,2}=\{\pi(i)=j, \pi(j)=i \},
\enas
and we may also write
\beas
&& A_1 = \{\pi(i)=i, \pi(j) \ne j \}, \ \ \ A_2 = \{\pi(i) \ne i, \pi(j) = j \} \\
&& A_3 = \{\pi(i) = j, \pi(j) \ne i \}, \ \ \ A_4 = \{\pi(j) = i, \pi (i) \ne j \},
\enas
and membership in $A_m$, $m=0,\ldots, 4$ depends only on $i,j,\pi(i),\pi(j)$.

Lastly, the sets $A_{5,m}$, $m=1,\ldots,4$ partition $A_5$, where
\beas
&& A_{5,1} = \{|i| = 2, |j| = 2, i \not\sim j \} \\
&& A_{5,2} = \{|i| = 2, |j| \ge 3 \}, \ \ \ A_{5,3} = \{|i| \ge 3, |j| = 2 \} \\
&& A_{5,4} = \{|i| \ge 3, |j| \ge 3 \} \cap A_5,
\enas
and membership in $A_{5,m}$, $m = 1,\ldots,4$ depends only on $i,j,\pi^{-1}(i),\pi^{-1}(j),\pi(i),\pi(j)$.
\end{lemma}

We now state Lemma 6.8 of \cite{CGS11}, noting that $y'$ and $y''$ there are incorrectly interchanged on the left hand side of \eqref{bdef}.
\begin{lemma}[\cite{CGS11}] \label{diff}
Let $\pi$ be a fixed permutation and $i \ne j$ distinct elements of $[n]$. Letting $\pi(-\alpha) =  \pi^{-1}(\alpha)$ for $\alpha \in [n]$ set
\beas
\chi_{i,j} = \{-i,-j,i,j \},  \text{ so that \ \ } \{\pi(\alpha), \alpha \in \chi_{i,j}\} = \{\pi^{-1}(i),\pi^{-1}(j),\pi(i),\pi(j) \}.
\enas
Then, for $\pi'' = \tau_{i,j} \pi' \tau_{i,j}$, and $y'$ and $y''$ given by \eqref{combdef} with $\pi'$ and $\pi''$ replacing $\pi$, respectively,
\bea \label{bdef}
y'-y'' = b(i,j,\pi(\alpha),\alpha \in \chi_{i,j})
\ena
where
\bea \label{sumdiff}
b(i,j,\pi(\alpha),\alpha \in \chi_{i,j}) = \sum_{m=0}^5 b_m(i,j,\pi(\alpha),\alpha \in \chi_{i,j}) \mathbf{1}_{A_m}
\ena
with $A_m$, $m=0,\ldots,5$ as in Lemma \ref{partition}, $b_0(i,j,\pi(\alpha),\alpha \in \chi_{i,j}) = 0$,
\bea \label{bidef}
&& \hspace{10pt} b_1(i,j,\pi(\alpha),\alpha \in \chi_{i,j}) = a_{i,i}+a_{\pi^{-1}(j),j}+a_{j,\pi(j)} -(a_{j,j}+a_{\pi^{-1}(j),i}+a_{i,\pi(j)}), \nn\\
&& \hspace{10pt} b_2(i,j,\pi(\alpha),\alpha \in \chi_{i,j}) = a_{j,j}+a_{\pi^{-1}(i),i}+a_{i,\pi(i)} -(a_{i,i}+a_{\pi^{-1}(i),j}+a_{j,\pi(i)}), \nn\\
&& \hspace{10pt} b_3(i,j,\pi(\alpha),\alpha \in \chi_{i,j}) = a_{\pi^{-1}(i),i}+a_{i,j}+a_{j,\pi(j)} -(a_{\pi^{-1}(i),j}+a_{j,i}+a_{i,\pi(j)}), \nn\\
&& \hspace{10pt} b_4(i,j,\pi(\alpha),\alpha \in \chi_{i,j}) = a_{\pi^{-1}(j),j}+a_{j,i}+a_{i,\pi(i)} -(a_{\pi^{-1}(j),i}+a_{i,j}+a_{j,\pi(i)}),\nn \\
&&\text{and} \nn\\
&& \hspace{10pt} b_5(i,j,\pi(\alpha),\alpha \in \chi_{i,j}) = a_{\pi^{-1}(i),i}+a_{i,\pi(i)}+a_{\pi^{-1}(j),j}+a_{j,\pi(j)} \nn\\
              && \hspace{140pt} -(a_{\pi^{-1}(i),j}+a_{j,\pi(i)}+a_{\pi^{-1}(j),i}+a_{i,\pi(j)}).
\ena
\end{lemma}

Next we construct an exchageable pair satisfying \eqref{approxsteinpair}, proved in the Appendix, which we call an approximate $\lambda,R$-Stein pair.

\begin{lemma} \label{combapproxstein}
For $n \ge 6$, let $\{a_{i,j}\}_{i,j=1}^n$ be an array of real numbers satisfying $a_{i,j} = a_{j,i}$ and $a_{\bullet,\bullet}=0$ where $a_{\bullet,\bullet}$ is as in \eqref{adotdot}. Let $\pi' \in \mathcal{S}_n$ a random permutation has the Ewens measure ${\cal E}_\theta$ with $\theta>0$. Further, let $I,J$ be chosen independently of $\pi$, uniformly from all pairs of distinct elements of $\{1,\ldots n\}$. Then, letting $\pi'' = \tau_{I,J}\pi' \tau_{I,J}$ and $Y'$ and $Y''$ be as in \eqref{combdef} with $\pi'$ and $\pi''$ replacing $\pi$, respectively, $(Y',Y'')$ is an approximate $4/n,R$-Stein pair with 
\bea \label{Rdef}
R(Y') =  \frac{1}{n(n-1)} \E [T|Y']
\ena
where
\bea \label{eq:defT}
T=
2(n+c_1-2(\theta+1)) \sum_{|i|=1}a_{i,i} + 2(c_1-2\theta)\sum_{|i| \ge 2} a_{i,i}-4\sum_{|i|,|j|=1, j \ne i}a_{i,j} -4\sum_{|i|=1,|j|\ge 2}a_{i,j}.
\ena
\end{lemma}

The next lemma provides bounds for $\E |R|$ and $|\E Y'R|$ that will be used when applying Theorems \ref{thm:approx.zb.dist.l1} and \ref{thm:approx.zb.dist.Kol}. To derive the bounds in Lemma \ref{remainderbounds}, we use consequences of Proposition \ref{c1moment}, which can be easily verified, that
\bea \label{c1moments}
&&\E \left[c_1\right] = \frac{\theta n }{\theta+n-1}, \ \ \ \ \E\left[ c_1(c_1-1)\right] = \frac{\theta^2 n_{(2)} }{(\theta+n-1)_{(2)}}, \nn  \\
&& \E \left[c_1^2\right] = \frac{\theta^2 n_{(2)} }{(\theta+n-1)_{(2)}}+\frac{\theta n }{\theta+n-1}, \text{ \ \  and } \nn\\
&& \E \left[c_1^2(c_1-1)^2\right] = \frac{\theta^4n_{(4)}}{(\theta+n-1)_{(4)}}+\frac{4\theta^3 n_{(3)}}{(\theta+n-1)_{(3)}}+\frac{2\theta^2n_{(2)}}{(\theta+n-1)_{(2)}}.
\ena

\begin{lemma} \label{remainderbounds}
Let $(Y',Y'')$ be an approximate $4/n,R$-Stein pair constructed as in Lemma \ref{combapproxstein} with $R$ as in \eqref{Rdef}. Then
\bea \label{Rbound}
\E|R| \le \frac{\theta M (12n+8\theta -10)}{(n-1)(\theta+n-1)} + \frac{2 \theta^2 M}{(\theta+n-1)_{(2)}}   ,
\ena 
and
\bea \label{YRbound}
|\E Y'R| \le \left(10\kappa_{\theta,n,1} + 4 \theta + \frac{4(\kappa_{\theta,n,1}(\theta+1)+\kappa_{\theta,n,2})}{n}  \right) \frac{ M \sigma}{n-1}
\ena
where $M = \max_{i,j} |a_{i,j} -a_{\bullet,\bullet}|$ with $a_{\bullet,\bullet}$ as in \eqref{adotdot} and $\kappa_{\theta,n,1}$ and $\kappa_{\theta,n,2}$ are given in \eqref{kappadef} and \eqref{kappadef2}, respectively.
\end{lemma}

\proof
By replacing $a_{i,j}$ by 
$a_{i,j} - a_{\bullet,\bullet}$ we may assume \eqref{zeroassumption} is satisfied, that is,
that $\E Y' =0$ and $a_{\bullet,\bullet} =0$ and thus demonstrate the claim with $M = \max_{i,j} |a_{i,j}|$.

We start with the first claim, using conditional Jensen's inequality and \eqref{Rdef} to obtain $|R| \le \E[|T|\,|Y']/(n(n-1))$ and therefore that $E|R| \le E|T|/(n(n-1))$, where $T$ is given by \eqref{eq:defT}. Using $|a_{i,j}| \le M$ for all $i,j$, we have
\beas
\E |R| &\le& \frac{2M \E |(n+c_1-2(\theta+1))c_1|}{n(n-1)} + \frac{2M\E|(c_1-2\theta)(n-c_1)|}{n(n-1)} \nn \\
               && \hspace{60pt} + \frac{4M\E[c_1(c_1-1)]}{n(n-1)} + \frac{4M\E[c_1(n-c_1)]}{n(n-1)} \nn \\
			 &=& \frac{2M \E |(n-1)c_1+c_1(c_1-1)-2\theta c_1|}{n(n-1)}  \nn \\
			   && \hspace{60pt}+ \frac{2M\E|(c_1-2\theta)(n-c_1)|}{n(n-1)}  + \frac{4M\E[(n-1)c_1]}{n(n-1)} \nn \\
			 &\le& \frac{2M \E [(n-1)c_1+c_1(c_1-1)+2\theta c_1]}{n(n-1)} + \frac{2M\E[nc_1 +2\theta n]}{n(n-1)} 
			     + \frac{4M\E[c_1]}{n} .
\enas

Now using \eqref{c1moments}, we obtain
\beas
\E|R| &\le& \frac{2\theta M }{\theta+n-1} + \frac{2\theta^2  M  }{(\theta+n-1)_{(2)}} +\frac{4\theta^2 M }{(n-1)(\theta+n-1)} \\
      && \hspace{10pt}+ \frac{2 \theta M n}{(n-1)(\theta+n-1)} + \frac{4 \theta M}{n-1}  + \frac{4 \theta M}{\theta+n-1}.
\enas
Combining similar terms in the last expression yields the bound in \eqref{Rbound}. 

Now we consider the second claim. First note that by \eqref{c1moments},
\beas
\E\left[\left(\sum_{|i|=1 } a_{i,i}\right)^2\right] \le M^2 \E[c_1^2] = \kappa_{\theta,n,1}^2 M^2 .
\enas
Hence, using that
\beas
n(n-1)|E[Y'R(Y')]|=|E[Y'E[T|Y']]|=|E[E[Y'T|Y']]|=|E[Y'T]|, 
\enas
for the first term in the product $Y'T$, bounding $c_1$ by $n$, we obtain

\bea \label{1stsecondterm}
\left|\frac{\E \left(2(n+c_1-2(\theta+1))Y' \sum_{|i|=1}a_{i,i}\right)}{n(n-1)}\right| &\le&
\left|\frac{4(n+\theta+1) \E \left(Y' \sum_{|i|=1 } a_{i,i}\right)}{n(n-1) }\right| \nn \\
&\le&  \frac{4\kappa_{\theta,n,1} M\sigma }{n-1}+\frac{4(\theta+1)\kappa_{\theta,n,1} M \sigma}{n(n-1)}
\ena
where here, and at similar steps later on, we apply the Cauchy--Schwarz inequality.

Now, to bound the second term in the $Y'T$ product, by \eqref{c1moments} we have
\beas
\E \left[c_1^2\left(\sum_{|i| \ge 2} a_{i,i}\right)^2\right] \le M^2n^2 \E[c_1^2] = \kappa_{\theta,n,1}^2 M^2 n^2 
\enas
and
\beas
\E \left[\left(\sum_{|i| \ge 2} a_{i,i}\right)^2\right] \le M^2n^2.
\enas
Thus, for that second term,
\bea \label{2ndsecondterm}
\left|\frac{\E\left((2c_1(\pi)-4\theta) Y' \sum_{|i| \ge 2} a_{i,i} \right)}{n(n-1)}\right| 
&\le& \left|\frac{2\E\left(c_1 Y' \sum_{|i| \ge 2} a_{i,i} \right)}{n(n-1)}\right| + \left|\frac{4\theta \E\left(Y' \sum_{|i| \ge 2} a_{i,i} \right)}{n(n-1)}\right| \nn \\
&\le& \frac{2 \kappa_{\theta,n,1} M\sigma }{n-1} + \frac{4\theta M \sigma}{n-1} .
\ena

For the third term, again by \eqref{c1moments},
\beas
\E \left[\left(\sum_{|i|=1,|j|=1, j \ne i}a_{i,j}\right)^2\right] \le M^2\E[c_1^2(c_1-1)^2] = M^2 \kappa_{\theta,n,2}^2,
\enas
and thus, 
\bea \label{3rdsecondterm}
\left|\frac{4\E\left(Y'\sum_{|i|=1,|j|=1, j \ne i}a_{i,j}\right)}{n(n-1)}\right| 
\le \frac{4\kappa_{\theta,n,2} M \sigma}{n(n-1)}.
\ena

For the fourth term, by \eqref{c1moments},
\beas
\E \left[\left(\sum_{|i|=1,|j|\ge 2}a_{i,j}\right)^2\right] \le n^2 M^2\E[c_1^2] = n^2 M^2 \kappa_{\theta,n,1}^2,
\enas
and thus, 
\bea \label{4thsecondterm}
\left|\frac{4\E\left(Y'\sum_{|i|=1,|j|\ge 2}a_{i,j}\right)}{n(n-1)}\right|
\le \frac{4 \kappa_{\theta,n,1} M \sigma}{n-1} .
\ena
Hence, summing \eqref{1stsecondterm} to \eqref{4thsecondterm}, we have

\beas
|\E Y'R| &\le& \frac{4\kappa_{\theta,n,1}M \sigma }{n-1}+\frac{4(\theta+1)\kappa_{\theta,n,1} M \sigma}{n(n-1)}
                                        + \frac{2 \kappa_{\theta,n,1} M \sigma}{n-1}\nn \\ 
																				&& \hspace{10pt}+ \frac{4\theta M \sigma}{n-1} 
																		    +\frac{4\kappa_{\theta,n,2} M \sigma}{n(n-1)}
																				+ \frac{4 \kappa_{\theta,n,1} M \sigma}{n-1}.
\enas
Combining similar terms and factoring $M\sigma/(n-1)$ out in the last expression yield the bound in \eqref{YRbound} and thus completes the proof.

\bbox

Next, Lemma \ref{combvar}, proved in the Appendix, provides detailed expressions for $\E[(Y'-Y'')^2]$ and $\sigma^2$. The order of $\sigma^2$ will be discussed in Remark \ref{sigmaorder}.

\begin{lemma} \label{combvar}
Let $n\ge 6$ and $(Y',Y'')$ be an approximate $4/n,R$-Stein pair constructed as in Lemma \ref{combapproxstein} with $R$ as in \eqref{Rdef}. Then 
\bea \label{Ydiff}
\E[(Y'-Y'')^2] = 2\beta_1+2\beta_3+\beta_{5,1}+2\beta_{5,2}+\beta_{5,4},
\ena
and
\bea \label{variance}
\sigma^2 = \frac{n}{8}(2\beta_1+2\beta_3+\beta_{5,1}+2\beta_{5,2}+\beta_{5,4})+\frac{n\E Y' R}{4},
\ena
where, with $b_m$, $m=1,\ldots,5$ are given in \eqref{bidef},
\bea \label{beta1}
\beta_1(\theta,n) = \frac{1}{n(n-1)(\theta+n-1)_{(3)}} \left(\theta^2 \sum_{i,j,s} b_1^2(i,j,i,s,i,s) + \theta \sum_{i,j,s,l} b_1^2(i,j,i,s,i,l)\right),
\ena
\bea \label{beta3}
\beta_3(\theta,n) = \frac{1}{n(n-1)(\theta+n-1)_{(3)}}  \left(\theta \sum_{i,j,r} b_3^2(i,j,r,i,j,r) + \sum_{i,j,r,l} b_3^2(i,j,r,i,j,l)\right),
\ena
\bea \label{beta51}
\beta_{5,1}(\theta,n) =  \frac{\theta^2}{n(n-1)(\theta+n-1)_{(4)}} \sum_{i,j,r,s} b_5^2(i,j,r,s,r,s),
\ena
\bea \label{beta52}
\beta_{5,2}(\theta,n) = \frac{\theta}{n(n-1)(\theta+n-1)_{(4)}} \sum_{i,j,r,s,l} b_5^2(i,j,r,s,r,l),
\ena
and
\bea \label{beta54}
\beta_{5,4}(\theta,n) &=& \frac{1}{n(n-1)(\theta+n-1)_{(4)}} \Bigg(\theta\sum_{i,j,r,k} b_5^2(i,j,r,k,k,r) \nn \\
                                         &&  \hspace{50pt}+\sum_{i,j,r,k,l} b_5^2(i,j,r,k,k,l) + \sum_{i,j,r,s,k,l} b_5^2(i,j,r,s,k,l)\Bigg).
\ena
\end{lemma}

\begin{remark} \label{sigmaorder}
Here we consider the order of $\sigma^2$ given in \eqref{variance}. Let $\beta_1$, $\beta_3$, $\beta_{5,1}$, $\beta_{5,2}$ and $\beta_{5,4}$ be given as in \eqref{beta1}-\eqref{beta54}, respectively. Using \eqref{YRbound}, we have
\beas
\sigma^2 &=& \frac{n}{8}(2\beta_1+2\beta_3+\beta_{5,1}+2\beta_{5,2}+\beta_{5,4})+\frac{n\E Y' R}{4}\\
         &\ge& \frac{n}{8}(2\beta_1+2\beta_3+\beta_{5,1}+2\beta_{5,2}+\beta_{5,4}) - \frac{n|\E Y' R|}{4} \\
				 &\ge& \frac{n}{8}(2\beta_1+2\beta_3+\beta_{5,1}+2\beta_{5,2}+\beta_{5,4}) \\
				&& \hspace{10pt}- \left(3\kappa_{\theta,n,1} + 1.2\theta + \frac{1.2(\kappa_{\theta,n,1}(\theta+1)+\kappa_{\theta,n,2})}{n}  \right) M \sigma. 
\enas
As discussed in Remark \ref{thetaremark}, $\kappa_{\theta,n,1}$ and $\kappa_{\theta,n,2}$ are of constant order in $n$ and thus $\sigma^2$ is of order $n$ whenever at least one of $\beta_\gamma(\theta,n)$, $\gamma = 1, 3, (5,1),(5,2),(5,4)$ are of constant order.
Since the final sum of \eqref{beta54} in the expression for $\beta_{5,4}(\theta,n)$ has $n_{(6)}$ terms and the denominator is of order $n^6$, $\beta_{5,4}(\theta,n)$ is of constant order in $n$ if the elements of $\{a_{i,j}: \{i,j\} \subset  [n] \}$ are chosen so that the values do not depend on $n$ and $b_5(i,j,r,s,k,l)$ with distinct $i,j,r,s,k,l \in [n]$ are nonzero for at least $\left\lfloor \delta n^6\right\rfloor$ terms for some $\delta>0$. 
For instance, if $a_{i,j}$, $\{i,j\} \subset [n]$ are independent identically distributed uniform random variables on $[0,1]$ then these sufficient conditions hold almost surely. 
\end{remark}

Now we are in the final step before proving Theorem \ref{L1theorem}, that is, to construct an approximate zero bias coupling that will be used when applying Theorems \ref{thm:approx.zb.dist.l1} and \ref{thm:approx.zb.dist.Kol}. Prior to doing that, we first specialize the outline in Lemma \ref{Wdagger2} to the more specific case where the random index ${\bf I}$ is chosen independently of the permutation and 
\bea \label{diffindex}
Y''-Y' = f({\bf I},\Xi_\alpha,\alpha \in \chi_{\bf I}),
\ena
where ${\bf I}$ and $\chi_{\bf I} \subset \chi$ are vectors of small dimensions and $f$ is a function with range being subset of $\mathbb{R}$, that is, we consider situations that $Y''-Y'$ depends on only a few variables.

Letting $Y',Y'',Y^\dagger,Y^\ddagger$ be constructed as in Lemma \ref{Wdagger2} satisfying \eqref{diffindex}, we follow Section 4.4.1 of \cite{CGS11} decomposing $P({\bf i},\xi_\alpha,\alpha \in \chi)$ as
\bea \label{pdecompose}
P({\bf i},\xi_\alpha,\alpha \in \chi) = P({\bf I} = {\bf i})P_{\bf i}(\xi_\alpha,\alpha \in \chi_{\bf i})P_{{\bf i}^c|{\bf i}}(\xi_\alpha,\alpha \not\in \chi_{\bf i}|\xi_\alpha,\alpha \in \chi_{\bf i}),
\ena
where $P_{\bf i}(\xi_{\alpha}, \alpha \in \chi_{\bf i})$ is the marginal distribution of $\xi_\alpha$ for $\alpha \in \chi_{\bf i}$, and $P_{{\bf i}^c|{\bf i}}(\xi_{\alpha}, \alpha \not\in \chi_{\bf i}|\xi_\alpha,\alpha \in \chi_{\bf i})$ the conditional distribution of $\xi_\alpha$ for $\alpha \not\in \chi_{\bf i}$ given $\xi_\alpha$ for $\alpha \in \chi_{\bf i}$.

For the square bias distribution, similarly, we decompose $P^\dagger({\bf i},\xi_{\alpha}, \alpha \in \chi)$ as
\bea \label{pdaggerdecompose}
P^\dagger({\bf i},\xi_\alpha,\alpha \in \chi) = P^\dagger({\bf I} = {\bf i})P^\dagger_{\bf i}(\xi_\alpha,\alpha \in \chi_{\bf i})P_{{\bf i}^c|{\bf i}}(\xi_\alpha,\alpha \not\in \chi_{\bf i}|\xi_\alpha,\alpha \in \chi_{\bf i})
\ena
where 
\bea \label{Idaggerdef}
P^\dagger({\bf I} = {\bf i}) = \frac{P({\bf I} = {\bf i}) \E f^2({\bf i},\Xi_\alpha,\alpha \in \chi_{\bf i})}{\E f^2({\bf I},\Xi_\alpha,\alpha \in \chi_{\bf I})},
\ena
and
\bea \label{dFdaggerdef}
P^\dagger_{\bf i}(\xi_{\alpha}, \alpha \in \chi_{\bf i}) = \frac{f^2({\bf i},\xi_\alpha,\alpha \in \chi_{\bf i})}{\E f^2({\bf i},\Xi_\alpha,\alpha \in \chi_{\bf i})}
                                                                 P_{\bf i}(\xi_{\alpha}, \alpha \in \chi_{\bf i}).
\ena

From this point, we denote ${\bf I}^\dagger$ for ${\bf I}$ that is generated from \eqref{Idaggerdef}. Notice that the representation of $P^\dagger_{\bf i}(\xi_{\alpha}, \alpha \in \chi)$ in \eqref{pdaggerdecompose} allows us to construct ${\bf I}^\dagger$ and $\{\Xi^\dagger_\alpha,\alpha \in \chi\}$ with distribution $P^\dagger({\bf i},\xi_\alpha,\alpha \in \chi)$ parallel to ${\bf I}$ and $\{\Xi_\alpha,\alpha \in \chi\}$ with distribution $P({\bf i},\xi_{\alpha}, \alpha \in \chi)$ in \eqref{pdecompose}. That is, we first choose ${\bf I}^\dagger$ by $P^\dagger({\bf I} = {\bf i})$ and then generate $\{\Xi^\dagger_\alpha,\alpha \in \chi_{\bf i}\}$ following distribution $P^\dagger_{\bf i}(\xi_{\alpha}, \alpha \in \chi_{\bf i})$ given ${\bf I}^\dagger = {\bf i}$. Finally, we generate $\{\Xi^\dagger_\alpha,\alpha \notin \chi_{\bf i}\}$ according to $P_{{\bf i}^c|{\bf i}}(\xi_{\alpha}, \alpha \not\in \chi_{\bf i}|\xi_\alpha,\alpha \in \chi_{\bf i})$. As the last term in \eqref{pdecompose} and \eqref{pdaggerdecompose} are exactly the same, the reader will see in the construction below that this equality will allow us to set $\Xi_\alpha = \Xi^\dagger_\alpha$ for most $\alpha$, yielding that $UY^\dagger+(1-U)Y^\ddagger$ and $Y'$ are close. 

\newpage

\begin{flushleft}
\textbf{Construction of an approximate zero bias coupling $(Y',Y^*)$:}
\end{flushleft}
Now we construct an approximate zero bias variable $Y^*$ from $Y'$, starting with its underlying permutations. Recall again that we consider the symbols $\pi$ and $Y$ interchageable with $\pi'$ and $Y'$, respectively. Let $\pi$ have the Ewens ${\cal E}_\theta$ distribution. By replacing $a_{i,j}$ by 
$a_{i,j} - a_{\bullet,\bullet}$ we may assume that \eqref{zeroassumption} is satisfied, that is,
that $\E Y' =0$ and $a_{\bullet,\bullet} =0$ and thus $M = \max_{i,j} |a_{i,j}|$.

Now we follow the outline described right before this construction to produce a coupling of $Y'$ to a pair $Y^\dagger,Y^\ddagger$, with the square bias distribution \eqref{eq:square.bias.pair} and then construct a coupling of $Y'$ to $Y^*$ satisfying \eqref{approxzerobias}, using uniform interpolation as in Lemma \ref{Wdagger}.

We first provide a brief overview of the construction, providing the full details later. To specialize the outline above to the case at hand, we let $f=b$ with $b$ given in Lemma \ref{diff}, ${\bf I} = (I,J)$, ${\bf i} = (i,j)$, $\Xi_\alpha = \pi(\alpha)$, $\Xi^\dagger_\alpha = \pi^\dagger(\alpha)$, $\chi = [n]$ and $\xi_\alpha, \alpha \in [n]$ denote distinct elements of $[n]$. As in Lemma \ref{diff}, we write $\pi(-\alpha) =  \pi^{-1}(\alpha)$ for $\alpha \in [n]$ and let $\chi_{\bf i} = \chi_{i,j} = \{-i,-j,i,j \}$ and
\beas
p_{i,j}(\xi_\alpha,\alpha \in \chi_{i,j}) = P(\pi(\alpha) = \xi_\alpha,\alpha \in \chi_{i,j}),
\enas
the distribution of the pre and post images of $i$ and $j$ under $\pi$.

By the decomposition in \eqref{pdecompose} and Lemma \ref{combapproxstein}, $\pi'$ and $\pi''$ can be constructed by choosing ${\bf I}$ with $P({\bf I} = {\bf i})$  uniformly on $[n]$, then constructing the pre and post images of $I$ and $J$ under $\pi'$ and the values of $\pi'$ on the remaining variables conditional on what has already been chosen and finally letting $\pi'' = \tau_{I,J}\pi'\tau_{I,J}$.

For the distribution of the pair $(Y^\dagger,Y^\ddagger)$ with the square bias distribution, we will first construct the underlying permutations $(\pi^\dagger,\pi^\ddagger)$. For this purpose, we follow the parallel decomposition of $P^\dagger({\bf i},\xi_\alpha,\alpha \in \chi)$ in \eqref{pdaggerdecompose} beginning with the indices ${\bf I}^\dagger = (I^\dagger,J^\dagger)$ with distribution \eqref{Idaggerdef}, 
\bea \label{IJgen}
P(I^\dagger=i,J^\dagger = j) = \frac{P(I=i,J=j)\E b^2(i,j,\pi(\alpha),\alpha \in \chi_{i,j})}{\E b^2(I,J,\pi(\alpha),\alpha \in \chi_{I,J})}
\ena 
where $b(i,j,\pi(\alpha),\alpha \in \chi_{i,j})$ is given in Lemma \ref{diff} and $\E b^2(I,J,\pi(\alpha),\alpha \in \chi_{I,J}) = \E[(Y'-Y'')^2]$ has been calculated in Lemma \ref{combvar}. Next, given $I^\dagger = i$ and $J^\dagger =  j$, we generate their pre and post images $\pi^{-\dagger}(I^\dagger),\pi^{-\dagger}(J^\dagger),\pi^{\dagger}(I^\dagger),\pi^{\dagger}(J^\dagger)$ with distribution \eqref{dFdaggerdef},
\bea \label{prepostprob}
p^\dagger_{i,j}(\xi_\alpha,\alpha \in \chi_{i,j}) = \frac{b^2(i,j,\xi_\alpha,\alpha \in \chi_{i,j})}{\E b^2(i,j,\pi(\alpha),\alpha \in \chi_{i,j})}p_{i,j}(\xi_\alpha,\alpha \in \chi_{i,j}).
\ena
Next we will construct the remaining images of $\pi^\dagger$ from $\pi$ conditional on what has already been chosen so that ${\bf I}^\dagger$ and $\pi^\dagger$ follow \eqref{pdaggerdecompose}, and that $\pi$ and $\pi^\dagger$ are close. Specializing the last factor of \eqref{pdaggerdecompose} to the case at hand, the remaining images of $\pi^\dagger$ given what has already been generated have distribution
\bea \label{remain}
p^\dagger_{{\bf i}^c|{\bf i}}\left(\xi_\alpha,\alpha \not\in \chi_{i,j}|\xi_\alpha,\alpha \in \chi_{i,j}\right) 
= P\left(\pi^\dagger(\alpha) = \xi_\alpha,\alpha \not\in \chi_{i,j}|\pi^\dagger(\alpha)=\xi_\alpha,\alpha \in \chi_{i,j}\right)
\ena
with $P$ follows the original Ewens $\mathcal{E}_\theta$ distribution. As the approximate $\lambda,R$-Stein pair $(Y',Y'')$ is defined as in \eqref{combdef} with $\pi$ replaced by $\pi'' = \tau_{I,J}\pi'\tau_{I,J}$ for $Y''$, we will then let $Y^\dagger$ and $Y^\ddagger$ be as in \eqref{combdef} with $\pi$ replaced by $\pi^\dagger$ and $\pi^\ddagger = \tau_{I^\dagger,J^\dagger}\pi^\dagger \tau_{I^\dagger,J^\dagger}$, respectively. Applying Lemma \ref{Wdagger2}, $(Y^\dagger,Y^\ddagger)$ will have the square bias distribution as in \eqref{eq:square.bias.pair}.

Next we construct $\pi^\dagger$ from the given $\pi'$. Recall from Section \ref{Ewens} that, with $B \subset [n]$, the reduced permutation $\sigma\setminus B$ acts on $[n]\setminus B$ and has cycle representation obtained by deleting all elements of $B$ in the cycle representation of $\sigma$. In view of Proposition \ref{conditionEwens}, we will keep the number of cycles $\#(\pi^\dagger \setminus \chi_{i,j})$ of the reduced permutation $\pi^\dagger \setminus \chi_{i,j}$, the same as $\#(\pi' \setminus \chi_{i,j})$ when applying \eqref{remain}.

For ease of notation, we denote the values
\bea  \label{prepost}
I^\dagger,J^\dagger,\pi^{-\dagger}(I^\dagger),\pi^{-\dagger}(J^\dagger),\pi^{\dagger}(I^\dagger),\pi^{\dagger}(J^\dagger)
\ena
generated in \eqref{IJgen} and \eqref{prepostprob} 
by $i,j,r,s,k,l$, respectively. By \eqref{IJgen}, $i \ne j$ and thus $r \ne s$ and $k \ne l$. Lemma \ref{partition} gives that, for $\pi^\dagger$, membership in $A_0,\ldots,A_4$ and $A_{5,1},\ldots,A_{5,4}$, defined in Lemma \ref{partition},  is determined by the generated $i,j,r,s,k,l$.

The remaining specification of $\pi^\dagger$ from $\pi$ depends on which case, or subcase, of the events $A_0,A_1$, $A_2$, $A_3$, $A_4$, $A_5$ is obtained. In each instance, we construct $\pi^\dagger$ from $\pi$ 
by removing $i,j,r,s$ from $\pi$'s cycle representation and inserting them into the resulting reduced permutation so that, respectively, $r$ and $k$ will be the pre and post images of $i$, and $s$ and $l$ will be those of $j$. As we only delete $i,j,r,s$ from $\pi$ and then insert these values without moving the remaining ones when constructing $\pi^\dagger$, it is clear that $\pi^\dagger \setminus \{i,j,r,s\}$ and $\pi \setminus \{i,j,r,s\}$ are identical. As the distribution of ${\bf I}^\dagger$ and their pre and post images follow \eqref{Idaggerdef} and \eqref{dFdaggerdef}, respectively, and the conditional distribution \eqref{remain} of the remaining values is the same as that for those of $\pi$, the pair $({\bf I}^\dagger,\pi^\dagger)$ follows the distribution in \eqref{pdaggerdecompose}.

\begin{figure}
  \includegraphics[width=\linewidth]{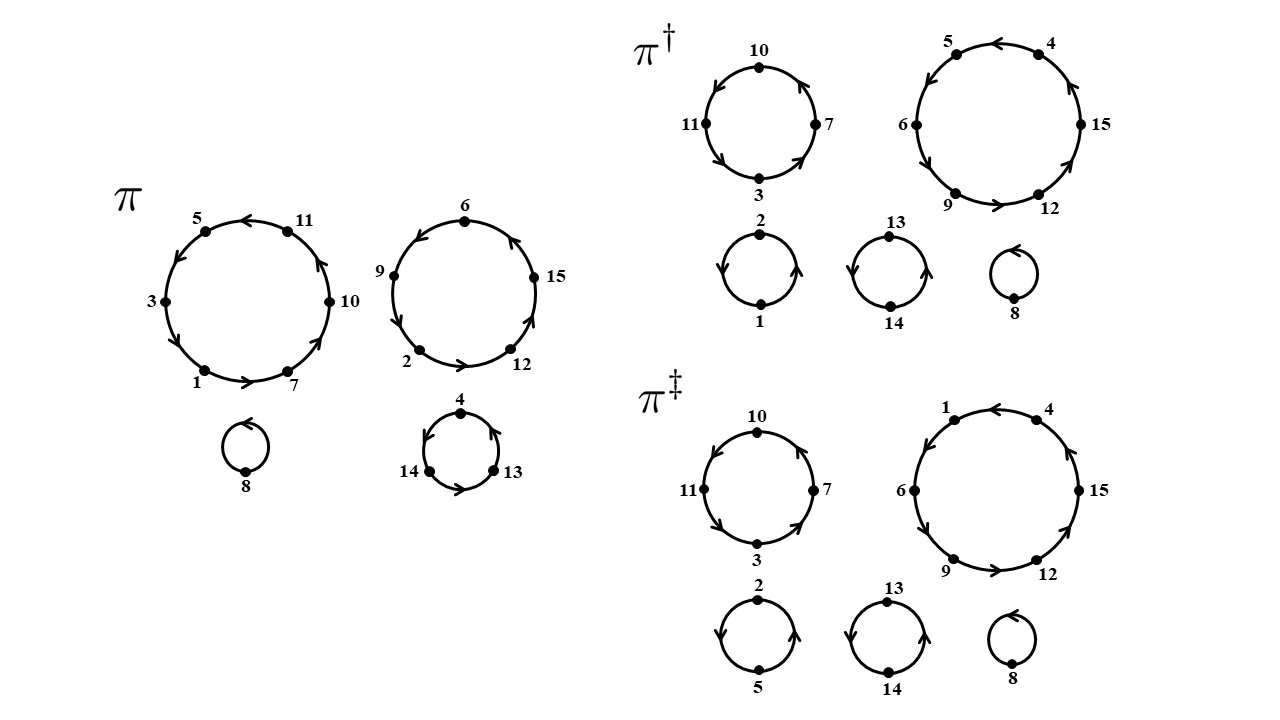}
  \caption{This figure shows an example of the construction of $\pi^\dagger$ and $\pi^\ddagger$ from $\pi$ when $n=15$ and $(i,j,r,s,k,l)=(1,5,2,4,2,6)$ that corresponds to the case $A_{5,2}$. One can easily notice that the reduced permutations $\pi\setminus \{1,5,2,4\}$, $\pi^\dagger \setminus \{1,5,2,4\}$, and $\pi^\ddagger \setminus \{1,5,2,4\}$ have the same cycle structure, that is, one 1-cycle, one 2-cycle and two 4-cycles.}
  \label{fig:example}
\end{figure}

In the following, we explicitly construct $\pi^\dagger$ from $\pi$ by separating into cases, or subcases, of the events $A_0,A_1$, $A_2$, $A_3$, $A_4$, $A_5$. To be more easily understandable, we also provide an example in Figure \ref{fig:example}. Although there are several cases, the idea in each case follows the same rule explained in the previous paragraph. Therefore, we only present the full detail in the first nonzero case $A_1$ and remark for the reader that the last case $A_{5,4}$ is the only one that contributes to the orders of the bounds in Theorem \ref{L1theorem} and hence it suffices for the reader to read only the first nonzero and the last cases. 

Below, we use the notation $i \rightarrow j$ if we change the cycle structure of $\pi$ so that $\pi^\dagger(i)=j$. We say `delete $i$ from $\pi$' if we delete $i$ from the cycle structure of $\pi$ and connect $\pi^{-1}(i) \rightarrow \pi(i)$ so that we end up with the reduced permutation $\pi \setminus \{i \}$. We say `insert $i$ in front of $k$' if we put $i$ between $\pi^{-1}(k)$ and $k$, that is, we end up with $\pi^{-1}(k) \rightarrow i \rightarrow k$ after the insertion. 

\textbf{Case $\mathbf{A_0}$}:
As $b_0 = 0$ by Lemma \ref{diff}, $p^\dagger_{i,j}$ as in \eqref{prepostprob} in the case $A_0$ is zero and therefore we need not consider this case. 

\textbf{Case $\mathbf{A_1}$}: We separate this case into two subcases, $A_{1,1} = \{|i|=1,|j|=2 \}$ and $A_{1,2} = \{|i|=1 , |j| \ge 3\}$. For $A_{1,1}$, we have $(i,j,r,s,k,l)=(i,j,i,s,i,s)$. We first recall that from \eqref{prepost} $r$ is the pre image of $i$ and thus $r=i$ means $i$ must be a fixed point. Similarly, since $s$ and $l$ are pre and post images of $j$, respectively, $l=s\ne j$ implies that $(s,j)$ must be a $2$-cycle. Hence, in this case, we simply delete $i,j,s$ from $\pi$, then let $i$ be a $1$-cycle and $(s,j)$ be a $2$-cycle. For $A_{1,2}$, we have $(i,j,r,s,k,l)=(i,j,i,s,i,l)$. By the same reasoning as for $A_{1,1}$, we keep $l$ at the original place and delete $i,j,s$ from $\pi$, then let $i$ be a $1$-cycle and insert $s \rightarrow j$ in front of $l$. For the remaining unmentioned values of the permutation, we keep them at the original places. Now we call the modified permutation $\pi^\dagger$. Since we have moved only $i,j$ and $s$, the reduced permutations $\pi^\dagger \setminus \{i,j,s \}$ and $\pi \setminus \{i,j,s \}$ are exactly the same. Applying Proposition \ref{conditionEwens} with $B=\{i,j,s \}$, we have that the remaining images of $\pi^\dagger$, conditioning on $\{\pi^\dagger(i)=i,\pi^\dagger(j)=l,\pi^\dagger(s)=j \}$ has distribution $p^\dagger_{{\bf i}^c|{\bf i}}(\xi_\alpha,\alpha \not\in \{-i,-j,i,j\}|\xi_{-i}=i,\xi_{-j}=s,\xi_i=i,\xi_j=l)$ in \eqref{remain}. Therefore the distribution of $({\bf I}^\dagger,\pi^\dagger)$ in $A_1$ case follows \eqref{pdaggerdecompose} and we have 
\beas
\pi(\phi) = \pi^\dagger(\phi) \text{ \ for all \ } \phi \not\in\mathcal{I}_{1} \text{ \ where \ } 
\mathcal{I}_1 = \{\pi^{-1}(\alpha),\beta:\alpha \in \{i,j,s,l\},\beta \in \{i,j,s\} \}. 
\enas
Hence $|\mathcal{I}_1|$ is at most $7$. From this point, for any index $\gamma$, we denote $\mathcal{I}_\gamma$ for the set that $\pi(\phi) = \pi^\dagger(\phi) $ for all $\phi \not\in\mathcal{I}_{\gamma}$ in the case $A_\gamma$. 

\textbf{Case $\mathbf{A_2}$}: Switching the role of $i$ and $j$, we follow the same construction as for $A_1$ and thus $|\mathcal{I}_2|$ is also at most $7$.

\textbf{Case $\mathbf{A_3}$}: Again we separate this case into $A_{3,1} = \{|i|=3,\pi^\dagger(i)=j \}$ and $A_{3,2} = \{|i|\ge 4 ,\pi^\dagger(i)=j\}$. For $A_{3,1}$, $(i,j,r,s,k,l)=(i,j,r,i,j,r)$, we delete $i,j,r$ from $\pi$ and let them be a $3$-cycle $r \rightarrow i \rightarrow j \rightarrow r$. For $A_{3,2}$, $(i,j,r,s,k,l)=(i,j,r,i,j,l)$, we delete $i,j,r$ from $\pi$ and insert $r\rightarrow i \rightarrow j$ in front of $l$. By the same reasoning as in $A_1$, the remaining images of $\pi^\dagger$, conditioning on $\{\pi^\dagger(i)=j,\pi^\dagger(j)=l,\pi^\dagger(r)=i \}$ has distribution $p^\dagger_{{\bf i}^c|{\bf i}}(\xi_\alpha,\alpha \not\in \{-i,-j,i,j\}|\xi_{-i}=r,\xi_{-j}=i,\xi_i=j,\xi_j=l)$ in \eqref{remain} and thus the distribution of $({\bf I}^\dagger,\pi^\dagger)$ in $A_3$ case follows \eqref{pdaggerdecompose}. In this case, 
\beas
\pi(\phi) = \pi^\dagger(\phi) \text{ \ for all \ } \phi \not\in\mathcal{I}_{3}  \text{ \ where \ } 
								\mathcal{I}_3 = \{\pi^{-1}(\alpha),\beta:\alpha \in \{i,j,r,l\},\beta \in \{i,j,r \} \}. 
\enas
Hence $|\mathcal{I}_3|$ is at most $7$. 

\textbf{Case $\mathbf{A_4}$}: Switching the role of $i$ and $j$, we follow the same construction as for $A_3$ and so $|\mathcal{I}_4|$ is at most $7$.

Moving on to $A_{5}$, here we separate it into $A_{5,1}$, $A_{5,2}$, $A_{5,3}$ and $A_{5,4}$, defined in Lemma \ref{partition}. 

\textbf{Case $\mathbf{A_{5,1}}$}: Here $(i,j,r,s,k,l)=(i,j,r,s,r,s)$, we delete $i,j,r,s$ from $\pi$ and then let both $(r,i)$ and $(s,j)$ be $2$-cycles. By the same argument as in the previous cases, the distribution of $({\bf I}^\dagger,\pi^\dagger)$ follows \eqref{pdaggerdecompose} and we end up with 
\beas
\pi(\phi) = \pi^\dagger(\phi) \text{ \ for all \ } \phi \not\in\mathcal{I}_{5,1}   \text{ \ where \ } \mathcal{I}_{5,1} = \{\pi^{-1}(\alpha),\alpha:\alpha \in \{i,j,r,s\} \}. 
\enas
Thus $|\mathcal{I}_{5,1}|$ is at most $8$.

\textbf{Case $\mathbf{A_{5,2}}$}: As $(i,j,r,s,k,l)=(i,j,r,s,r,l)$, we delete $i,j,r,s$ from $\pi$ and then let $(r,i)$ be a $2$-cycle and insert $s\rightarrow j$ in front of $l$. In this case, again the distribution of $({\bf I}^\dagger,\pi^\dagger)$ follows \eqref{pdaggerdecompose} and 
\beas
\pi(\phi) = \pi^\dagger(\phi) \text{ \ for all \ } \phi \not\in\mathcal{I}_{5,2} 
\text{ \ where \ } \mathcal{I}_{5,2} = \{\pi^{-1}(\alpha),\beta:\alpha \in \{i,j,r,s,l\},\beta \in \{i,j,r,s \} \}. 
\enas
Hence $|\mathcal{I}_{5,2}|$ is at most $9$. 

\textbf{Case $\mathbf{A_{5,3}}$}: Switching the role of $i$ and $j$, we follow the same construction as for $A_{5,2}$ and thus $|\mathcal{I}_{5,3}|$ is at most $9$.

\textbf{Case $\mathbf{A_{5,4}}$}: We again need to break this case into subcases depending on the size of $\{i,j,r,s,k,l \}$. As the case that contributes the most to the difference between $\pi$ and $\pi^\dagger$ is the case where $i,j,r,s,k,l$ are all distinct, we only discuss the construction of this particular case. We delete $i,j,r,s$ from $\pi$ and then insert $r \rightarrow i$ and $s \rightarrow j$ in front of $k$ and $l$, respectively. In this case, the distribution of $({\bf I}^\dagger,\pi^\dagger)$ again follows \eqref{pdaggerdecompose} and 
\beas
\pi(\phi) = \pi^\dagger(\phi) \text{ \ for all \ } \phi \not\in\mathcal{I}_{5,4} 
 \text{ \ where \ } \mathcal{I}_{5,4} = \{\pi^{-1}(\alpha),\beta:\alpha \in \{i,j,r,s,k,l\}, \beta \in \{i,j,r,s\}\}. 
\enas
Thus $|\mathcal{I}_{5,4}|$ is at most $10$.

Now the construction of $\pi^\dagger$ has been specified in every case and subcase and the distribution of $({\bf I}^\dagger,\pi^\dagger)$ follows \eqref{pdaggerdecompose}. Therefore, setting
\beas
\pi^\ddagger = \tau_{I^\dagger,J^\dagger} \pi^\dagger \tau_{I^\dagger,J^\dagger} 
\enas
results in a collection of variables $I^\dagger,J^\dagger$ and a pair of permutations with the square bias distribution \eqref{squarebias}. Hence, letting $Y',Y^\dagger,Y^\ddagger$ given by \eqref{combdef} with $\pi'$, $\pi^\dagger$ and $\pi^\ddagger$, respectively, yields a coupling of $Y'$ to the variables $Y^\dagger,Y^\ddagger$ with the square bias distribution as required in Lemma \ref{Wdagger}. 
Invoking Lemma \ref{Wdagger} to $Y^*=UY^\dagger+(1-U)Y^\ddagger$ with $U \sim \mathcal{U}[0,1]$ be independent of $Y^\dagger,Y^\ddagger$, we have $Y^*$ be an approximate $Y'$-zero bias variable satisfying \eqref{approxzerobias}, as desired. 

The following lemma provides a bound of the difference between $Y^*$ and $Y'$.
\begin{lemma} \label{Ystarbound}
Let $Y'$ and $M$ be defined as in the statement of Theorem \ref{L1theorem} and $Y^*$ be constructed as in the construction above. Then
$|Y^* - Y'| \le 20 M$.
\end{lemma}
\proof
By the same argument as in the proof of Lemma 6.10 in \cite{CGS11}, we have
\beas
|Y^* - Y'| \le 2 \max_{\gamma \in \Gamma}|\mathcal{I}_\gamma| M,
\enas
where $\Gamma = \{1,2,3,4,(5,1),(5,2),(5,3),(5,4)\}$. Then the claim follows from that $|\mathcal{I}_\gamma|$ is at most $10$ from the construction above.

\bbox

Now we have all ingredients to prove our main theorem.

\bigskip

\noindent {\bf Proof of Theorem \ref{L1theorem}} 
As before, by replacing $a_{i,j}$ by 
$a_{i,j} - a_{\bullet,\bullet}$ we may assume that \eqref{zeroassumption} is satisfied, that is,
that $\E Y' =0$ and $a_{\bullet,\bullet} =0$ and thus proceed the construction with $M = \max_{i,j} |a_{i,j}|$. First we construct an approximate $4/n,R$-Stein pair $(Y',Y'')$ as in Lemma \ref{combapproxstein} with the remainder $R$ given in \eqref{Rdef}. Then we construct an approximate zero bias variable $Y^*$ satisfying \eqref{approxzerobias} as in the construction right above this proof.

Now we apply Theorems \ref{thm:approx.zb.dist.l1} and \ref{thm:approx.zb.dist.Kol}, handling three terms on the right hand side of \eqref{l1approx} and \eqref{Kolapprox}. We note that the three terms from the two theorems are different only on their constants so we handle both $L^1$ and $L^\infty$ upper bounds at the same time. For the first term, by Lemma \ref{Ystarbound}, we have
\bea \label{firstmain}
\frac{1}{\sigma}|Y^*-Y'| \le \frac{20M}{\sigma}.
\ena

Now we handle the last two terms containing the remainder $R$. Applying \eqref{YRbound} and \eqref{Rbound} in Lemma \ref{remainderbounds}, respectively, and using that $\lambda = 4/n$ and $n \ge 6$, we obtain
\bea \label{secondmain}
\frac{|\E Y'R|}{\sigma^2 \lambda} \le \left(3\kappa_{\theta,n,1} + 1.2 \theta + \frac{\kappa_{\theta,n,1}(\theta+1)+\kappa_{\theta,n,2}}{n-1} \right) \frac{M}{\sigma} ,
\ena
and
\bea \label{lastmain}
\frac{\E|R|}{\sigma \lambda} \le 
   \frac{0.6M\theta (6n+4\theta-5)}{\sigma(\theta+n-1)} + \frac{0.5\theta^2Mn}{\sigma(\theta+n-1)_{(2)}},
\ena
where $\kappa_{\theta,n,1}$ and $\kappa_{\theta,n,2}$ are given in \eqref{kappadef} and \eqref{kappadef2}, respectively.

Invoking Theorems \ref{thm:approx.zb.dist.l1} and \ref{thm:approx.zb.dist.Kol}, using \eqref{firstmain}, \eqref{secondmain} and \eqref{lastmain}, now yields the $L^1$ and $L^\infty$ upper bounds in \eqref{mainbound} and \eqref{mainbound2}, respectively.

Next we follow the idea in \cite{Eng81} to prove the $L^\infty$ lower bound in \eqref{mainlowerbound} in the case that $a_{i,j}$ are all integers. By Chebyshev's inequality,
\beas
P(\E Y - \sqrt{3}\sigma < Y < \E Y + \sqrt{3}\sigma) \ge 2/3.
\enas
The random variable $Y$ is discrete having at most $n!$ possibilities and for each pair of possibilities, the difference between them is at least $1$ as $a_{i,j}$ are all integers.
Therefore the interval $(\E Y - \sqrt{3}\sigma,\E Y + \sqrt{3}\sigma)$ contains at most $\left\lceil (2\sqrt{3}\sigma+1) \right\rceil$ possible values of $Y$ which implies that the largest point mass $p(n)$ in the distribution of $Y$ satisfies
\beas
p(n) \ge \frac{2}{6\sqrt{3}\sigma +3}.
\enas
As $\Phi(x)$ the distribution function of $Z$ is continuous, we have
\beas
\sup_{\-\infty < x < \infty} |P(Y<x) - \Phi((x-\E Y)/\sigma)| \ge \frac{1}{2} p(n) = \frac{1}{6\sqrt{3} \sigma+3}.
\enas

\bbox

\section{Appendix} \label{appendix}

We prove Lemmas \ref{combapproxstein} and \ref{combvar} in this section.

\noindent {\bf Proof of Lemma \ref{combapproxstein}}
As Ewens measure is constant on cycle type, the exchangeability claim follows immediately from the proof of Lemma 6.9 of \cite{CGS11}. 

It remains to show that $Y',Y''$ satisfies \eqref{approxsteinpair} with $R$ given by \eqref{Rdef}. 
As $Y'$ is a function of $\pi'$ the tower property of conditional expectation yields that $$\E[Y'-Y''|Y'] = \E[\E[Y'-Y''| \pi]|Y'],$$ and we begin by computing the conditional expectation given $\pi$ of the difference \\$\sum_{m=0}^5 b_m(i,j,\pi(\alpha),\alpha \in \chi_{i,j}) \mathbf{1}_{A_m}$ in \eqref{sumdiff} with $A_0,\ldots,A_5$ given in Lemma \ref{partition}, with $i,j$ replaced by $I,J$. 

First we have that $b_0 = 0$. As given in (6.85) of \cite{CGS11}, the contribution to $n(n-1)\E[Y'-Y''| \pi]$ from $b_1$ and $b_2$ totals to
\bea \label{6.85}
2(n-c_1)\sum_{|i|=1} a_{i,i} + 4c_1 \sum_{|i|\ge 2} a_{i,\pi(i)} 
      - 2c_1\sum_{|i|\ge 2} a_{i,i} -2 \sum_{|i|=1,|j|\ge 2} a_{i,j} -2 \sum_{|i|\ge 2,|j|=1}a_{i,j},
\ena
and likewise (6.88) shows $b_3$ and $b_4$ contribute
\bea \label{6.88}
6\sum_{|i|\ge 3}a_{i,\pi(i)} - 4\sum_{|i|\ge 3} a_{\pi^{-1}(i),\pi(i)} -2\sum_{|i| \ge 3} a_{\pi(i),i},
\ena
and (6.90) shows the first four terms of $b_5$ contribute
\bea \label{6.90}
4(n-2-c_1)\sum_{|i|=2}a_{i,\pi(i)} + 4(n-3-c_1) \sum_{|i| \ge 3} a_{i,\pi(i)}.
\ena
Lastly, the contribution from the fifth term of $b_5$ is given by (6.91), and separating out the cases where $i\not =j$ and $i=j$ in the first sum there, that expression can be seen equivalent to 
\bea \label{6.91}
-\sum_{|i|\ge 4}\sum_{j \sim i,j \ne i}a_{i,j}+ \sum_{|i| \ge 4}(a_{i,\pi(i)}+a_{\pi^{-1}(i),\pi(i)}) - \sum_{|i|\ge 2,|j|\ge 2}\sum_{j\not\sim i}a_{i,j}.
\ena

To simplify \eqref{6.91}, let $a \wedge b = \min(a,b)$, and follow (6.92) of \cite{CGS11}, separating out the cases where $j=i$ and $j \not =i$, resulting in the identity
\bea \label{6.92}
\theta \sum_{i=1}^n a_{i,i} + \sum_{i \ne j} a_{i,j} = \theta\sum_{|i|\ge 1} a_{i,i}  + \sum_{|i|\ge 4}\sum_{j \sim i,j \ne i}a_{i,j} + \sum_{|i|\le 3}\sum_{j \sim i,j \ne i}a_{i,j} \nn\\
+ \sum_{|i|\ge 2,|j|\ge 2}\sum_{j\not\sim i}a_{i,j} + \sum_{|i|\wedge |j|=1} \sum_{j \not\sim i}a_{i,j}.
\ena 
Since $\theta\sum_{i=1}^n a_{i,i} + \sum_{i \ne j} a_{i,j} = n(\theta+n-1)a_{\bullet,\bullet} = 0$ by the assumption, we replace the sum of the first and last terms in \eqref{6.91} by the sum of  the first, third and fifth terms in \eqref{6.92}. Now \eqref{6.91} equals 
\beas
\theta\sum_{|i|\ge 1} a_{i,i} + \sum_{|i|\le 3}\sum_{j \sim i,j \ne i}a_{i,j} + \sum_{|i| \ge 4}(a_{i,\pi(i)}+a_{\pi^{-1}(i),\pi(i)}) + \sum_{|i|\wedge |j|=1} \sum_{j \not\sim i}a_{i,j}.
\enas
Using that $\pi^2(i) = \pi^{-1}(i)$ when $|i|=3$ and there are only two points in the cycle when $|i|=2$, we obtain
\beas
\theta \sum_{|i|\ge 1}a_{i,i} + \sum_{|i| \ge 2}a_{i,\pi(i)}+ \sum_{|i| \ge 3}a_{\pi^{-1}(i),\pi(i)} + \sum_{|i|\wedge |j|=1} \sum_{j \not\sim i}a_{i,j}.
\enas
Combining this contribution with the next three terms of $b_5$, each of which yields the same amount, gives the total
\bea \label{6.93}
4\theta \sum_{|i|\ge 1}a_{i,i} + 4\sum_{|i| \ge 2}a_{i,\pi(i)}+ 4\sum_{|i| \ge 3}a_{\pi^{-1}(i),\pi(i)} + 4\sum_{|i|\wedge |j|=1} \sum_{j \not\sim i}a_{i,j}.
\ena

Combining \eqref{6.93} with the contribution \eqref{6.90} from the first four terms in $b_5$, the $b_1$ and $b_2$ terms in \eqref{6.85} and the $b_3$ and $b_4$ terms \eqref{6.88}, yields
\beas
 && n(n-1)\E[Y'-Y''|\pi] \nn \\
 && \hspace{20pt}             = 4(n-1) \sum_{|i|=2} a_{i,\pi(i)} + (4n-2) \sum_{|i|\ge 3} a_{i,\pi(i)}  -2 \sum_{|i|\ge 3} a_{\pi(i),i} \nn \\
  && \hspace{40pt}                   +2(n-c_1+2\theta) \sum_{|i|=1} a_{i,i} - 2(c_1-2\theta)\sum_{|i| \ge 2} a_{i,i} \nn\\
	 && \hspace{40pt}								  +4\sum_{|i|\wedge |j|=1, j \not\sim i}a_{i,j} -2\sum_{|i|=1,|j|\ge 2}a_{i,j}-2\sum_{|i|\ge 2,|j|= 1}a_{i,j} \nn \\
    && \hspace{20pt}  =4(n-1) \sum_{|i|=2} a_{i,\pi(i)} + (4n-2) \sum_{|i|\ge 3} a_{i,\pi(i)}  -2 \sum_{|i|\ge 3} a_{\pi(i),i} +4(n-1) \sum_{|i|=1} a_{i,i} \nn\\  
 && \hspace{40pt}	                    - 2(n+c_1-2(\theta+1)) \sum_{|i|=1}a_{i,i} - 2(c_1-2\theta)\sum_{|i| \ge 2} a_{i,i} \nn\\
	&& \hspace{40pt}									  +4\sum_{|i|\wedge |j|=1, j \not\sim i}a_{i,j} -2\sum_{|i|=1,|j|\ge 2}a_{i,j}-2\sum_{|i|\ge 2,|j|= 1}a_{i,j},
\enas
where the final two expressions are identical but for a rewriting of the coefficient of $\sum_{|i|=1}a_{i,i}$.

Now applying the assumption that $a_{i,j}=a_{j,i}$ to the third term, using that $\pi(i) = i$ for $|i|=1$ for the fourth term, and $a_{i,j}=a_{j,i}$ again to obtain
\beas
&&4\sum_{|i|\wedge |j|=1, j \not\sim i}a_{i,j} -2\sum_{|i|=1,|j|\ge 2}a_{i,j}-2\sum_{|i|\ge 2,|j|= 1}a_{i,j}\\
&&\hspace{50pt}= 4\sum_{|i|=1,|j|=1, j \ne i}a_{i,j} +4\sum_{|i|=1,|j|\ge 2}a_{i,j}
\enas
we have
\beas
n(n-1)\E[Y'-Y''|\pi] 
&=& 4(n-1) \sum_{|i|=2} a_{i,\pi(i)} + (4n-2) \sum_{|i|\ge 3} a_{i,\pi(i)}  \\
&& \hspace{2pt} - 2 \sum_{|i|\ge 3} a_{i,\pi(i)} +4(n-1) \sum_{|i|=1} a_{i,\pi(i)}-T\\  
&=& 4(n-1)Y' - T,
\enas
where $T$ is given in \eqref{eq:defT}. Now conditioning on $Y'$ proves the claim that $\E[Y''|Y'] = (1-\lambda) Y' + R(Y') $ with $\lambda = 4/n$ and $R(Y')$ as in \eqref{Rdef}.

\bbox

\bigskip

In the following proof, for ease of notation, we write $\Sigma_{\alpha_1,\ldots,\alpha_m}$ for the sum over distinct $\alpha_k \in [n]$.

\noindent {\bf Proof of Lemma \ref{combvar}}
We first calculate $\E[(Y'-Y'')^2]$ with the help of Lemma \ref{diff}. Note that we write $b_m(i,j,\pi(\alpha),\alpha \in \chi_{i,j}) = b_m(i,j,\pi^{-1}(i),\pi^{-1}(j),\pi(i),\pi(j))$ in what follows. As $b_0=0$, moving on to $A_1$, and recalling that the summations in this proof are over distinct values, we have 
\beas
\E[(Y'-Y'')^2]\mathbf{1}_{A_1} &=& \E b_1^2(I,J,\pi^{-1}(I),\pi^{-1}(J),\pi(I),\pi(J)) \mathbf{1}_{\{\pi(I)=I,\pi(J) \ne J\}} \\
&=& \frac{1}{n(n-1)(\theta+n-1)_{(3)}} \left(\theta^2 \sum_{i,j,s} b_1^2(i,j,i,s,i,s) + \theta \sum_{i,j,s,l} b_1^2(i,j,i,s,i,l)\right),
\enas 
noting that there are $n(n-1)$ possibilities for $I$ and $J$, applying Proposition \ref{probEwens} with $B=\{i,j,s \}$ we see that
the factor $\theta^2/(\theta+n-1)_{(3)}$ in the first term is the probability that $\pi(i)=i$, $\pi(j)=s$ and $\pi(s)=j$ and the factor $\theta/(\theta+n-1)_{(3)}$ in the second term is the probability that $\pi(i)=i$, $\pi(j)=l$ and $\pi(s) = j$.  By symmetry, $A_2$ contributes the same amount. 

For $A_3$, by similar reasoning, we have,
\beas
\E[(Y'-Y'')^2]\mathbf{1}_{A_3} &=& \E b_3^2(I,J,\pi^{-1}(I),\pi^{-1}(J),\pi(I),\pi(J)) \mathbf{1}_{\{\pi(I)=J,\pi(J) \ne I\}} \\
&=& \frac{1}{n(n-1)(\theta+n-1)_{(3)}}  \left(\theta \sum_{i,j,r} b_3^2(i,j,r,i,j,r) + \sum_{i,j,r,l} b_3^2(i,j,r,i,j,l)\right).
\enas 
The contribution from $A_4$ is the same as that from $A_3$.

Lastly, consider the contribution from $A_5$. Starting with $A_{5,1}$, we have
\beas
\E[(Y'-Y'')^2]\mathbf{1}_{A_{5,1}} &=& \E b_5^2(I,J,\pi^{-1}(I),\pi^{-1}(J),\pi(I),\pi(J)) \mathbf{1}_{\{|I|=2,|J|=2,I \not\sim J\}} \\
&=& \frac{\theta^2}{n(n-1)(\theta+n-1)_{(4)}} \sum_{i,j,r,s} b_5^2(i,j,r,s,r,s).
\enas 

For $A_{5,2}$, we have 
\beas
\E[(Y'-Y'')^2]\mathbf{1}_{A_{5,2}} &=& \E b_5^2(I,J,\pi^{-1}(I),\pi^{-1}(J),\pi(I),\pi(J)) \mathbf{1}_{\{|I|=2,|J|\ge 3\}} \\
&=& \frac{\theta}{n(n-1)(\theta+n-1)_{(4)}} \sum_{i,j,r,s,l} b_5^2(i,j,r,s,r,l).
\enas 
Again by symmetry, $A_{5,3}$ contributes the same.

For $A_{5,4}$, we have
\beas
\E[(Y'-Y'')^2]\mathbf{1}_{A_{5,4}} &=& \E b_5^2(I,J,\pi^{-1}(I),\pi^{-1}(J),\pi(I),\pi(J)) \mathbf{1}_{\{|I|\ge 3,|J|\ge 3\}} \\
&=& \frac{1}{n(n-1)(\theta+n-1)_{(4)}} 
\Bigg(\theta  \sum_{i,j,r,k} b_5^2(i,j,r,k,k,r) \\
&&\hspace{10pt}+\sum_{i,j,r,k,l} b_5^2(i,j,r,k,k,l) + \sum_{i,j,r,s,k,l} b_5^2(i,j,r,s,k,l)\Bigg).
\enas 
Summing all cases, yields $\E[(Y'-Y'')^2]$ in \eqref{Ydiff}.

Since $\E[(Y'-Y'')^2] = 2\left(\lambda \sigma^2 -\E Y'R \right)$ and $\lambda = 4/n$, we have
\beas
\sigma^2 &=& \frac{n\E[(Y'-Y'')^2]}{8} + \frac{n\E Y' R}{4} \\
         &=& \frac{n}{8}(2\beta_1+2\beta_3+\beta_{5,1}+2\beta_{5,2}+\beta_{5,4})+\frac{n\E Y' R}{4},
\enas
where $\beta_1$, $\beta_3$, $\beta_{5,1}$, $\beta_{5,2}$ and $\beta_{5,4}$ are defined as in \eqref{beta1}, \eqref{beta3}, \eqref{beta51}, \eqref{beta52} and \eqref{beta54}, respectively.

\bbox

\bigskip

\section*{Acknowledgements} 
This work is part of my Ph.D. thesis at the University of Southern California. I thank my advisor, Larry Goldstein, for his invaluable guidance and fruitful discussions upon finishing this research. I also thank Jason Fulman for suggesting this interesting topic. Finally, I thank an anonymous referee for pointing out an issue in the previous version.

\bigskip

\def\cprime{$'$}

\end{document}